\documentclass[11pt]{article}
\input{psfig.sty}
\usepackage{amsbsy}
\usepackage{amsmath}
\usepackage{amssymb}

\renewcommand{\baselinestretch} {1.3}
\makeatletter 
\def\singlespace{\def\baselinestretch{1}\@normalsize}

\@addtoreset{equation}{section}
\renewcommand{\theequation} {\arabic{section}.\arabic{equation}}

\newtheorem{definition}{Definition}

\newtheorem{lemma}{Lemma}
\newtheorem{proposition}{Proposition}

\def\cal{\mathcal}
\def\hat{\widehat}
\def\tilde{\widetilde}

\def\bbeta{{\boldsymbol{\beta}}}

\def\btheta{{\boldsymbol{\theta}}}
\def\bevector{\boldsymbol{e}}
\def\calS{\cal S}
\def\bX{{\bf X}}
\def\bx{{\bf x}}
\def\bzero{{\bf 0}}
\def\ID{\mathrm{I}}
\def\bvectorX{{\texttt{X}}}
\def\bvectorx{{\texttt{x}}}
\def\bvectorZ{{\texttt{Z}}}
\def\bvectorz{{\texttt{z}}}
\def\bA{{\boldsymbol{A}}}

\def\bz{\mbox{\bf z}}
\def\br{\mbox{\bf r}}
\def\bB{{\bf B}}
\def\bH{{\bf H}}

\def\bW{\mbox{\bf W}}
\def\bK{\mbox{\bf K}}
\def\bV{\mbox{\bf V}}
\def\LB{\mathrm{LB}}

\def\bu{{\mathbf u}}
\def\calH{\cal H}
\def\calK{\cal K}
\def\calL{\cal L}
\def\calX{\cal X}
\def\calC{\cal C}
\def\ta{\texttt{a}}
\def\logit{\mathrm{logit}}
\def\tr{\mathrm{tr}}
\def\diag{\mathrm{diag}}
\def\var{{\mathrm{var}}}
\def\cov{{\mathrm{cov}}}
\def\MISE{{\mathrm{MISE}}}
\def\MPEC{{\mathrm{MPEC}}}
\def\hAMISE{h_{\mathrm{AMISE}}}
\def\hAMPEC{h_{\mathrm{AMPEC}}}
\def\hmin{h_{\mathrm{min}}}

\def\CV{\mathrm{CV}}
\def\ACV{\mathrm{ACV}}
\def\ECV{\mathrm{ECV}}
\def\hACV{\hat h_{\ACV}}
\def\hECV{\hat h_{\ECV}}
\def\etal{{et al.}\ }
\def\Err{{\mathrm{Err}}}
\def\err{{\mathrm{err}}}
\def\bI{{\bf I}}
\def\sign{\mathrm{sign}}

\begin{document}

\title
{\bf Assessing Prediction Error of Nonparametric Regression and
Classification under Bregman Divergence}
\author
{
Jianqing Fan \\
Department of Operations Research \\
and Financial Engineering \\
Princeton University \\
Princeton, NJ 08544 \\
jqfan@princeton.edu \and
Chunming Zhang \\
Department of Statistics \\
University of Wisconsin \\
Madison, WI 53706-1685 \\
cmzhang@stat.wisc.edu}
\date{\today}

\maketitle
%\begin{singlespace}
%\begin{footnotetext}
%{}
%\end{footnotetext}
%\end{singlespace}
\begin{abstract}

Prediction error is critical to assessing the performance of
statistical methods and selecting statistical models.  We propose
the cross-validation and approximated cross-validation methods for
estimating prediction error under a broad $q$-class of Bregman
divergence for error measures which embeds nearly all of the
commonly used loss functions in regression, classification
procedures and machine learning literature.  The approximated
cross-validation formulas are analytically derived, which facilitate
fast estimation of prediction error under the Bregman divergence. We
then study a data-driven optimal bandwidth selector for the
local-likelihood estimation that minimizes the overall prediction
error or equivalently the covariance penalty. It is shown that the
covariance penalty and cross-validation methods converge to the same
mean-prediction-error-criterion. We also propose a lower-bound
scheme for computing the local logistic regression estimates and
demonstrate that it is as simple and stable as the local
least-squares regression estimation. The algorithm monotonically
enhances the target local-likelihood and converges. The idea and
methods are extended to the generalized varying-coefficient models
and semiparametric models.
\end{abstract}


\noindent {\bf Short title}: Assessing Prediction Error under
Bregman Divergence

\newpage

\section{Introduction}

Assessing prediction error lies at the heart of statistical model
selection and forecasting.  Depending on the needs of statistical
learning and model selection, the quadratic loss is not always
appropriate. In binary classification, for example, the
misclassification error rate is more suitable.  The corresponding
loss, however, does not differentiate the predictive powers of two
procedures which forecast the class label being $1$ with
probabilities $60\%$ and $90\%$ respectively.  When the true class
label is $1$, the second procedure is more accurate, whereas when
the true class label is $0$, the first procedure is more accurate.
The quantification of the predicative accuracy requires an
appropriate introduction of loss functions. An example of this is
the negative Bernoulli log-likelihood loss function. Other important
margin-based loss functions have been introduced for binary
classification in the machine learning literature (Hastie,
Tibshirani and Friedman 2001). Hence, it is important to assess the
prediction error under a broader class of loss functions.

A broad and important class of loss functions is the Bregman
$q$-class divergence.  It accounts for different types of output
variables and includes the quadratic loss, the deviance loss for
exponential family of distributions, misclassification loss and
other popular loss functions in machine learning.  See Section
\ref{sect-2}. Once a prediction error criterion is chosen, the
estimates of prediction error are needed.  Desirable features
include computational expediency and theoretical consistency. In the
traditional nonparametric regression models, the residual-based
cross-validation (CV) is a useful data-driven method for the
automatic smoothing (Wong 1983; Rice 1984; H\"ardle, Hall, and
Marron 1992; Hall and Johnstone 1992) and can be handily computed.
With the arrival of the optimism theorem (Efron 2004), estimating
the prediction error becomes estimating covariance-penalty terms.
Following Efron (2004), the covariance-penalty can be estimated
using model-based bootstrap procedures.  A viable model-free method
is the cross-validated estimation of the covariance-penalty. Both
methods can be shown to be asymptotically equivalent to the first
order approximation. However, both methods are extremely
computationally intensive in the context of local-likelihood
estimation, particularly for the large sample sizes.  The challenge
then arises from efficient computation of the estimated prediction
error based on the cross-validation.

The computational problem is resolved via the newly developed
approximate formulas for the cross-validated covariance-penalty
estimates. A key component is to establish the ``leave-one-out
formulas" which offer an analytic connection between the
leave-one-out estimates and their ``keep-all-in'' counterparts. This
technical work integrates the infinitesimal perturbation idea
(Pregibon 1981) with the Sherman-Morrison-Woodbury formula (Golub
and Van Loan 1996, p. 50).  It is a natural extension of the
cross-validation formula for least-squares regression estimates, and
is applicable to both parametric and nonparametric models.

The applications of estimated prediction error pervade almost every
facet of statistical model selection and forecasting. To be more
specific, we focus on the local-likelihood estimation in varying
coefficient models for response variables having distributions in
the exponential family. Typical examples include fitting the
Bernoulli distributed binary responses, and the Poisson distributed
count responses, among many other non-normal outcomes. As a flexible
nonparametric model-fitting technique, the local-likelihood method
possesses nice sampling properties. For details, see, for example,
Tibshirani and Hastie (1987), Staniswalis (1989), Severini and
Staniswalis (1994), and Fan, Heckman, and Wand (1995). An important
issue in application is the choice of smoothing parameter.
Currently, most of the existing methods deal with the Gaussian type
of responses; clearly there is a lack of methodology for
non-Gaussian responses. The approximate cross-validation provides a
simple and fast method for this purpose. The versatility of the
choice of smoothing parameters is enhanced by an appropriate choice
of the divergence measure in the $q$-class of loss functions.

The computational cost of the approximate CV method is further
reduced via a newly introduced empirical version of CV, called ECV,
which is based on an empirical construction of the ``degrees of
freedom", a notion that provides useful insights into the
local-likelihood modeling complexity. We propose a data-driven
bandwidth selection method, based on minimizing ECV, which will be
shown to be asymptotically optimal in minimizing a broad $q$-class
of prediction error. Compared with the two-stage bandwidth selector
of Fan, Farmen, and Gijbels (1998), our proposed method has a
broader domain of applications  and can be more easily understood
and implemented.

Some specific attentions are needed for the local logistic
regression with binary responses, whose distribution belongs to an
important member of the exponential family. To address the numerical
instability, we propose to replace the Hessian matrix by its global
lower-bound (LB) matrix, which does not involve estimating parameter
vectors and therefore can easily be inverted before the start of the
Newton-Raphson (NR) iteration. A similar idea of LB was used in
B\"ohning and Lindsay (1988) for some parametric fitting. We make a
conscientious effort to further develop this idea for the local
logistic estimation. The resulting LB method gains a number of
advantages: The LB algorithm, at each iteration, updates the
gradient vector but does not recalculate the Hessian matrix, thus is
as simple and stable as the local least-squares regression
estimation. The LB method ensures that each iterative estimate
monotonically increases the target local-likelihood. In contrast,
this property is not shared by the standard NR method. Hence, the LB
iteration is guaranteed to converge to the true local MLE, whereas
the NR is not necessarily convergent. Moreover, we develop a new and
adaptive data-driven method for bandwidth selection which can
effectively guard against under- or over-smoothing.

The paper is organized as follows. Section \ref{sect-2} addresses
the issue of estimating prediction error. Sections \ref{sect-3}
develops computationally feasible versions of the cross-validated
estimates of the prediction error. Section \ref{sect-4} proposes a
new bandwidth selection method for binary responses, based on the LB
method and the cross-validated estimates of the prediction error.
Sections \ref{sect-5}--\ref{sect-6} extend our results to
generalized varying-coefficient models and semiparametric models
respectively. Section \ref{sect-7} presents simulation evaluations
and Section \ref{sect-8} analyzes real data. Technical conditions
and proofs are relegated to the Appendix.

\section{Estimating Prediction Error} \label{sect-2}

To begin with, we consider that the response variable $Y$ given the
vector $\bvectorX$ of input variables has a distribution in the
exponential family, taking the form,
\begin{eqnarray} \label{b1}
f_{Y|\bvectorX}(y; \theta(\bvectorx))=\exp[\{y
\theta(\bvectorx)-b(\theta(\bvectorx))\}/a(\psi)+c(y, \psi)],
\end{eqnarray}
for some known functions $a(\cdot)$, $b(\cdot)$, and $c(\cdot,
\cdot)$, where $\theta(\bvectorx)$ is called a canonical parameter
and $\psi$ is called a dispersion parameter, respectively. It is
well known that
\begin{eqnarray*}
m(\bvectorx) \equiv
E(Y|\bvectorX=\bvectorx)=b'(\theta(\bvectorx)),\quad
\mbox{and}\quad \sigma^2(\bvectorx) \equiv
\var(Y|\bvectorX=\bvectorx)=a(\psi) b''(\theta(\bvectorx)).
\end{eqnarray*}
See Nelder and Wedderburn (1972) and McCullagh and Nelder (1989).
The canonical link is $g(\cdot)=(b')^{-1}(\cdot)$, resulting in
$g(m(\bvectorx))=\theta(\bvectorx)$.  For simplicity of notation
and exposition, we will focus only on estimating the canonical
parameter $\theta(\bvectorx)$.  The results can easily be
generalized to other link functions.

\subsection{Bregman Divergence} \label{sect-2.1}

\begin{figure}[htbp]
\centerline{\psfig{figure=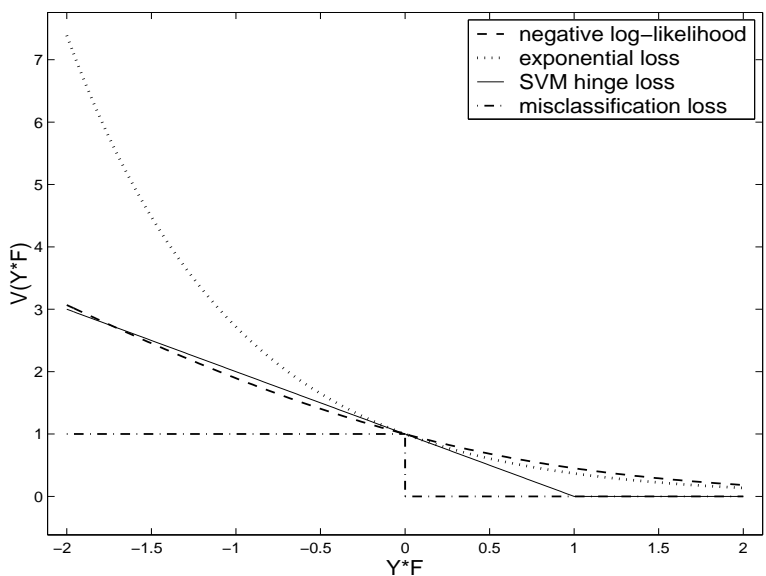,height=3in,width=4in}}
\begin{singlespace}
\caption {\textsl{Illustration of Margin-Based Loss Functions. Line
types are indicated in the legend box. Each function has been
re-scaled to pass the point $(0,1)$.}} \label{figure-1}
\end{singlespace}
\end{figure}

The prediction error depends on the divergence measure.  For
non-Gaussian responses, the quadratic loss function is not always
adequate.  For binary classification, a reasonable choice of
divergence measure is the misclassification loss, $Q(Y, \hat m) =
\ID\{Y \not = \ID(\hat m > .5)\}$, where $\ID(\cdot)$ is an
indicator function and $\hat m$ is an estimator. However, this
measure does not differentiate the predictions $\hat m = .6$ and
$\hat m = .9$ when $Y = 1$ or $0$. In the case that $Y = 1$, $\hat m
= .9$ gives a better prediction than $\hat m = .6$. The negative
Bernoulli log-likelihood, $Q(Y, \hat m) = - Y \ln (\hat m) - (1-Y)
\ln (1-\hat m)$, captures this. Other loss functions possessing
similar properties include the hinge loss function, $Q(Y, \hat m) =
\max\{1 - (2Y-1)\sign(\hat m-.5), 0 \}$, in the support vector
machine and the exponential loss function, $Q(Y, \hat m) = \exp\{ -
(Y-.5) \ln(\hat m/(1-\hat m))\}$, popularly used in {AdaBoost}.
These four loss functions, shown in Figure \ref{figure-1}, belong to
the margin-based loss functions written in the form, $V(Y^*F)$, for
$Y^*=2Y-1$ and some function $F$.
\begin{center}
{[ \textsl{Put Figure \ref{figure-1} about here} ]}
\end{center}

\begin{figure}[htbp]
\centerline{\psfig{figure=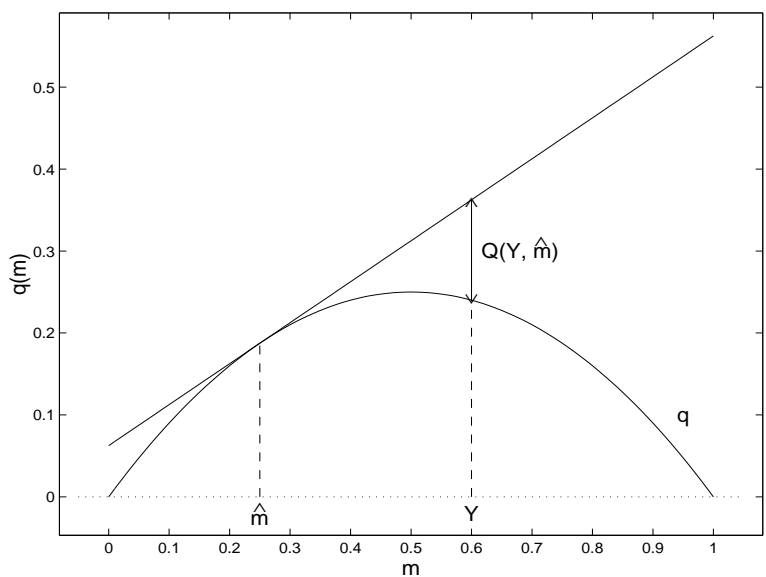,height=3in,width=4in}}
\begin{singlespace}
\caption {\textsl{Illustration of the Bregman Divergence $Q(Y,\hat
m)$ as Defined in (\ref{b2}). The concave curve is the $q$-function;
the two dashed lines give locations of $Y$ and $\hat m$; the solid
strict line is $q(\hat m)+q'(\hat m)(Y-\hat m)$; the vertical line
with arrows at each end is $Q(Y,\hat m)$.}} \label{figure-2}
\end{singlespace}
\end{figure}

To address the versatility of loss functions, we appeal to a device
introduced by Bregman (1967). For a concave function $q(\cdot)$,
define a $q$-class of error measures $Q$ as
\begin{eqnarray} \label{b2}
    Q(Y,\hat m)=q(\hat m)+q'(\hat m)(Y-\hat m)-q(Y).
\end{eqnarray}
A graphical illustration of $Q$ associated with $q$ is displayed in
Figure~\ref{figure-2}. Due to the concavity of $q$, $Q$ is
non-negative. However, since $Q(\cdot,\cdot)$ is not generally
symmetric in its arguments, $Q$ is not a ``metric" or ``distance" in
the strict sense. Hence, we call $Q$ the Bregman ``divergence" (BD).
\begin{center}
{[ \textsl{Put Figure \ref{figure-2} about here} ]}
\end{center}

It is easy to see that, with the flexible choice of $q$, the BD is
suitable for a broad class of error measures. Below we present some
notable examples of the $Q$-loss constructed from the $q$-function.
\begin{itemize}
\item A function $q_1(m) = a m - m^2$ for some constant $a$  yields
the quadratic loss $Q_1(Y, \hat {m}) = (Y - \hat {m})^2$.

\item For the exponential family $(\ref{b1})$,
the function $q_2(m) = 2 \{ b(\theta) - m \theta \}$ with
$b'(\theta) = m$ results in  the deviance loss,
\begin{eqnarray} \label{b3}
   Q_2 (Y, \hat {m}) = 2 \{Y (\tilde{\theta} -
\hat {\theta}) - b(\tilde{\theta}) + b(\hat {\theta}) \},
\end{eqnarray}
where $b'(\tilde{\theta}) = Y$ and $b'(\hat {\theta}) =  \hat m$.

\item For a binary response variable $Y$, the function $q(m) =
\min\{m, (1-m)\}$ gives the misclassification loss; the function
$q(m) = .5 \min\{m, (1-m)\}$ results in the hinge loss; the function
$q_3(m) = 2\{m (1-m)\}^{1/2}$ yields the exponential loss,
\begin{eqnarray} \label{b4}
Q_3 (Y, \hat m) = \exp \{ -(Y-.5) \ln (\hat m/(1-\hat m)) \}.
\end{eqnarray}

\end{itemize}

\subsection{Prediction Error under Bregman Divergence}

Let $m_i=m(\bvectorX_i)$ and $\hat m_i$ be its estimate based on
independent observations $\{(\bvectorX_i, Y_i)\}_{i=1}^{n}$. Set
\begin{eqnarray*}
\err_i = Q(Y_i,\hat m_i)
\quad \mbox{and} \quad
\Err_i = E_o\{ Q(Y_i^o,\hat m_i)\},
\end{eqnarray*}
where $Y_i^o$ is an independent copy of $Y_i$ and is independent
of $(Y_1,\ldots,Y_n)$, and $E_o$ refers to the expectation with
respect to the probability law of $Y_i^o$.  Note that the
conditional prediction error, defined by $\Err_i$, is not
observable, whereas the apparent error, $\err_i$, is observable.
As noted in Tibshirani (1996), directly estimating the
{conditional} prediction error is very difficult. Alternatively,
estimating $\Err_i$ is equivalent to estimating the difference
$O_i=\Err_i-\err_i$, called \textsl{optimism}.

Efron (2004) derives the optimism theorem to represent the
expected optimism as the covariance penalty, namely, $E(O_i) =
2\,\cov(\hat \lambda_i,Y_i)$, where $\hat \lambda_i = -q'(\hat
m_i)/2$.
 As a result, the predictive error can be
estimated by
\begin{eqnarray} \label{b5}
\hat \Err_i=\err_i+2\,\hat \cov_i,
\end{eqnarray}
where $\hat \cov_i$ is an estimator of the covariance penalty,
$\cov(\hat \lambda_i,Y_i)$. This is an insightful generalization of
AIC.  Henceforth, the total prediction error $\Err=\sum_{i=1}^n
\Err_i$ can be estimated by $\hat \Err=\sum_{i=1}^n \hat \Err_i$.

\subsection{Estimation of Covariance Penalty}

In nonparametric estimation, we write $\hat m_{i,h}$ and $\hat
\lambda_{i,h}$ to stress their dependence on a smoothing parameter
$h$. According to (\ref{b5}), the total prediction error is given by
\begin{eqnarray} \label{b6}
\hat \Err (h)= \sum_{i=1}^{n} Q(Y_i,\hat m_{i,h})+\sum_{i=1}^n
2\,\hat \cov(\hat \lambda_{i,h},Y_i).
\end{eqnarray}
Estimating the covariance-penalty terms in (\ref{b6}) is not
trivial. Three approaches are potentially applicable to the
estimation of the covariance penalty: model-based bootstrap
developed in Efron (2004), the data-perturbation (DP) method
proposed by Shen, Huang and Ye (2004), and model-free
cross-validation.  All three are computationally intensive in the
context of local-likelihood estimation (the DP method also needs to
select a perturbation size). In contrast, the third method allows us
to develop approximate formulas to significantly gain computational
expedience. For this reason, we focus on the cross-validation
method.

The cross-validated estimation of  $\Err_i$ is $Q(Y_i,\hat m_{i,
h}^{-i})$, where the superscript $-i$ indicates the deletion of the
$i$th data point $(\bvectorX_i,Y_i)$ in the fitting process. This
yields the cross-validated estimate of the total prediction error by
\begin{eqnarray} \label{b7}
{\hat \Err}^{\mathrm{CV}}(h)= \sum_{i=1}^n Q(Y_i,\hat
m_{i,h}^{-i}).
\end{eqnarray}
Naive computation of $\{ \hat m_{i,h}^{-i} \}$ is intensive. Section
\ref{sect-3} will devise strategies by which actual computations of
the leave-one-out estimates are not needed.  A distinguished feature
is that our method is widely applicable to virtually all regression
and classification problems.  The approximated CV is particularly
attractive to a wide array of large and complex problems in which a
quick and crude selection of the model parameter is needed.

By comparing (\ref{b7}) with (\ref{b6}), the covariance penalty in
(\ref{b7}) is estimated by
$$
   \sum_{i=1}^n \{Q(Y_i,\hat m_{i,h}^{-i}) - Q(Y_i,\hat m_{i,h})\}.
$$
This can be linked with the jackknife method for estimating the
covariance penalty.  Hence, it is expected that the
cross-validation method is asymptotically equivalent to a
bootstrap method.

\subsection{Asymptotic Prediction Error} \label{sect-2.4}

To gain insight on $\hat \Err(h)$, we appeal to the asymptotic
theory. Simple algebra shows that $\Err_i(h) = Q(m_i,\hat
m_{i,h})+E\{Q(Y_i,m_i)\}$. By Taylor's expansion and (\ref{b2}),
$$
 Q(m_i,\hat m_{i,h})\doteq -(\hat m_{i,h}-m_i)^2q''(\hat
 m_{i,h})/2.
$$
Hence,
$$
\Err_i(h) \doteq -(\hat m_{i,h}-m_i)^2q''(\hat
m_{i,h})/2+E\{Q(Y_i,m_i)\}.
$$
Note that the last term does not depend on $h$ and hence $\hat
\Err(h)$ is asymptotically equivalent to the
mean-prediction-error-criterion,
\begin{eqnarray} \label{b8}
\MPEC(h)= - 2^{-1}\int E[\{\hat m_h(x)-m(x)\}^2|\calX] q''(m(x))
f_X(x) dx,
\end{eqnarray}
with $\calX=(\bvectorX_1,\ldots,\bvectorX_n)$ and $f_X(x)$ being the
probability density of $\bvectorX$. This criterion differs from the
mean-integrated-squared-error criterion defined by
\begin{eqnarray} \label{b9}
\MISE(h)= \int E[\{\hat m_h(x)-m(x)\}^2|\calX]
\{b''(\theta(x))\}^{-2} f_X(x) dx,
\end{eqnarray}
recalling that
$$
  \hat{\theta}(x) - \theta(x)  \doteq \{b''(\theta(x))\}^{-1} \{\hat
  m_h(x)-m(x)\}.
$$

Expression (\ref{b8}) reveals that asymptotically, different loss
functions automatically introduce different weighting schemes in
(\ref{b8}). This provides a useful insight into various error
measures used in practice.  The weighting schemes vary substantially
over the choices of $q$. In particular, for the $q_1$-function
yielding the quadratic-loss in Section \ref{sect-2.1}, we have
\begin{eqnarray*}
\MPEC_1(h)=\int E[\{\hat m_h(x)-m(x)\}^2|\calX] f_X(x) dx.
\end{eqnarray*}
For the $q_2$-function producing the deviance-loss, we have
\begin{eqnarray*}
\MPEC_2(h) = \int E[\{\hat m_h(x)-m(x)\}^2|\calX]
\{b''(\theta(x))\}^{-1} f_X(x) dx.
\end{eqnarray*}
For the $q_3$-function inducing the exponential-loss for the binary
responses, we have
\begin{eqnarray*}
\MPEC_3(h) = \int E[\{\hat m_h(x)-m(x)\}^2|\calX] \frac{f_X(x)}{4
[m(x)\{1-m(x)\}]^{3/2}} dx.
\end{eqnarray*}

\section{Approximate Cross-Validation} \label{sect-3}

This section aims at deriving the approximate and empirical versions
of (\ref{b7}) for the local maximum likelihood estimator. We focus
on the univariate problem in this section.  The results will be
extended to the generalized varying coefficient models in Section
\ref{sect-5} incorporating multivariate covariates.

Assume that the function $\theta(\cdot)$ has a $(p+1)$-th
continuous derivative at a point $x$.  For $X_j$ close to $x$, the
Taylor expansion implies that
\begin{eqnarray*}
\theta(X_j)\doteq  \bx_j(x)^T\bbeta(x),
\end{eqnarray*}
in which $\bx_j(x)=(1,(X_j-x), \ldots, (X_j-x)^p)^T$ and
$\bbeta(x)=(\beta_0(x), \ldots, \beta_p(x))^T$. Based on the
independent observations, the local parameters can be estimated by
maximizing the local log-likelihood,
\begin{eqnarray} \label{c1}
\ell(\bbeta; x) \equiv \sum_{j=1}^n l(\bx_j(x)^T\bbeta; Y_j)
K_h(X_j-x),
\end{eqnarray}
in which $l(\cdot;y)=\ln \{f_{Y|X}(y;\cdot)\}$ denotes the
conditional log-likelihood function, $K_h(\cdot)=K(\cdot/h)/h$ for
a kernel function $K$ and $h$ is a bandwidth. Let $\hat
\bbeta(x)=(\hat \beta_0(x), \ldots, \hat \beta_p(x))^T$ be the
local maximum likelihood estimator. Then, the local MLEs of
$\theta(x)$ and $m(x)$ are given by $\hat{\theta}(x) =
\hat{\beta}_0(x)$ and $\hat m(x)=b'(\hat \theta(x))$,
respectively. A similar estimation procedure, based on the $n-1$
observations excluding $(X_i,Y_i)$, leads to the local
log-likelihood function $\ell^{-i}(\bbeta; x)$, and the
corresponding local MLEs, $\hat \bbeta^{-i}(x)$, $\hat
\theta^{-i}(x)$, and $\hat m^{-i}(x)$,  respectively.

\subsection{Weighted Local-Likelihood} \label{sect-3.1}

To compute approximately $\hat \bbeta^{-i}(x)$ from $\hat
\bbeta(x)$, we apply the ``infinitesimal perturbation" idea
developed in Pregibon (1981). We introduce the weighted local
log-likelihood function,
\begin{eqnarray} \label{c2}
\ell_{i,\delta}(\bbeta; x)=\sum_{j=1}^n \delta_{ij}\,
l(\bx_j(x)^T\bbeta; Y_j) K_h(X_j-x),
\end{eqnarray}
with the weight $\delta_{ii}= \delta$ and the rest weights
$\delta_{ij}=1$. Let $\hat \bbeta_{i, \delta}(x)$ be the maximizer.
Note that when $\delta = 1$, this estimator is the local maximum
likelihood estimator and when $\delta =0$, it is the leave-one-out
estimator.

The weighted local MLE is usually found via the Newton-Raphson
iteration,
\begin{eqnarray} \label{c3}
\bbeta_{L}=\bbeta_{L-1}- \{\triangledown^2
\ell_{i,\delta}(\bbeta_{L-1}; x)\}^{-1} \triangledown
\ell_{i,\delta}(\bbeta_{L-1}; x), \quad L=1, 2, \ldots,
\end{eqnarray}
where $\triangledown \ell$ denotes the gradient vector and
$\triangledown^2 \ell$ the Hessian matrix. (Explicit expressions of
$\triangledown \ell$ and $\triangledown^2 \ell$ are given in Lemma
\ref{lem-2}.) When the initial estimator $\bbeta_0$ is good enough,
the one-step ($L=1$) estimator is as efficient as the fully iterated
estimator (Fan and Chen 1999).

The key ingredient for calculating the leave-one-out estimator is
to approximate it by its one-step estimator using the
``keep-all-in'' estimator $\hat{\bbeta}(x)$ as the initial value.
Namely, $\hat{\bbeta}^{-i}(x)$ is approximated by (\ref{c3}) with
$\bbeta_0 = \hat{\bbeta}(x)$ and $\delta=0$.  With this idea and
other explicit formulas and approximation, we can derive the
approximate leave-one-out formulas.

\subsection{Leave-One-Out Formulas} \label{sect-3.2}

 Let $\bX(x)=(\bx_1(x),\ldots$,$\bx_n(x))^T$,
\begin{eqnarray} \label{c4}
\bW(x; \bbeta) = \diag\{K_h(X_j-x) b''(\bx_j(x)^T\bbeta)\},
\end{eqnarray}
and $S_n(x; \bbeta)=\bX(x)^T \bW(x;\bbeta)\bX(x)$.  Define
\begin{eqnarray} \label{c5}
\calH(x;\bbeta) = \{\bW(x;\bbeta)\}^{1/2} \bX(x)
\{S_n(x;\bbeta)\}^{-1} \bX(x)^T \{\bW(x;\bbeta)\}^{1/2}.
\end{eqnarray}
This projection matrix is an extension of the hat matrix in the
multiple regression and will be useful for computing the
leave-one-out estimator.  Let $\calH_{ii}(x;\bbeta)$ be its
$i$-th diagonal element and $H_{i} =\calH_{ii}(X_i;\hat
\bbeta(X_i))$.  Then, our main results can be summarized as
follows.
\begin{proposition} \label{prop1}
Assume condition $(A2)$ in the Appendix. Then for $i=1, \ldots, n$,
\begin{eqnarray}
\hat \bbeta^{-i}(x)-\hat \bbeta(x) &\doteq& -\frac {\{S_n(x;\hat
\bbeta(x))\}^{-1} \bx_i(x) K_h(X_i-x) \{Y_i-b'(\bx_i(x)^T \hat
\bbeta(x))\}}
         {1-\calH_{ii}(x;\hat \bbeta(x))}, \label{c6} \\
\hat \theta^{-i}_i-\hat \theta_i &\doteq& -\frac {H_i} {1-H_i}
         (Y_i - \hat m_i)/b''(\hat \theta_i), \label{c7} \\
\hat m^{-i}_i-\hat m_i &\doteq& -\frac {H_i}{1 - H_i} (Y_i-\hat
m_i). \label{c8}
\end{eqnarray}
\end{proposition}

Note that the approximation becomes exact when the loss function is
the quadratic loss.  In fact, Zhang (2003) shows an explicit
delete-one-out formula,
$$
    \hat m^{-i}_i = \hat m_i- \frac{H_i}{1 - H_i}  (Y_i-\hat m_i).
$$
In addition, (\ref{c6})--(\ref{c8}) hold for $h \to \infty$, namely,
the parametric maximum likelihood estimator.  Even the results for
this specific case appear new.  Furthermore, the results can easily
be extended to the estimator that minimizes the local Bregman
divergence, replacing $l(\bx_j(x)^T \bbeta; Y_j)$ in (\ref{c1}) by
$Q(Y_j, g^{-1}(\bx_j(x)^T \bbeta))$.

 Using  Proposition~\ref{prop1}, we can derive a simplified formula for
computing the cross-validated estimate of the overall prediction
error.

\begin{proposition} \label{prop2}
Assume conditions $(A1)$ and $(A2)$ in the Appendix. Then
\begin{eqnarray}
{\hat \Err}^{\CV} \doteq \sum_{i=1}^n \Big[Q(Y_i,\hat m_{i})
+2^{-1}q''(\hat m_i)(Y_i-\hat m_i)^2
\big\{1-{1}/{(1-H_i)^{2}}\big\} \Big]. \label{c9}
\end{eqnarray}
\end{proposition}

Proposition \ref{prop2} gives an approximation formula, which avoids
computing ``leaving-one-out'' estimates, for all $q$-class of loss
functions. In particular, for the function $q_1$, we have
$$
    \sum_{i=1}^{n} (Y_i-\hat m_i^{-i})^2 \doteq \sum_{i=1}^{n} (Y_i-\hat
    m_i)^2/(1- H_i)^2.
$$
For this particular loss function,  the approximation is actually
exact. For the function $q_2$ leading to the deviance loss $Q_2$
defined in (\ref{b3}), we have
\begin{eqnarray}
\sum_{i=1}^n Q_2(Y_i,\hat m_i^{-i}) \doteq \sum_{i=1}^n
\bigg[Q_2(Y_i,\hat m_i)- \frac{(Y_i-\hat m_i)^2}{b''(\hat
       \theta_i)} \big\{1-{1}/{(1-H_i)^{2}}\big\}\bigg]. \label{c10}
\end{eqnarray}
For the exponential loss defined in (\ref{b4}) for binary
classification, we observe
\begin{eqnarray*}
\sum_{i=1}^n Q_3(Y_i,\hat m_i^{-i}) \doteq \sum_{i=1}^n
\bigg[Q_3(Y_i,\hat m_i)- \frac{(Y_i-\hat m_i)^2}{4
\{\hat{m}_i(1-\hat{m}_i)\}^{3/2}}
\big\{1-{1}/{(1-H_i)^{2}}\big\}\bigg].
\end{eqnarray*}

\subsection{Two Theoretical Issues} \label{sect-3.3}

Two theoretical issues are particularly interesting. The first one
concerns the asymptotic convergence of $\hACV$, the minimizer of the
right hand side of (\ref{c9}). Following a suitable modification to
the result of Altman and MacGibbon (1998), the ratio $\hACV/\hAMPEC$
converges in probability to $1$, where $\hAMPEC$ is the minimizer of
the asymptotic form of $\MPEC(h)$ defined in (\ref{b8}).

The explicit expression of $\hAMPEC$, associated with the $q$-class
of error measures, can be obtained by the delta method. Setting
$-2^{-1}q''(m(x))\{b''(\theta(x))\}^2$ to be the weight function,
$\hAMPEC$ (for odd degrees $p$ of local polynomial fitting) can be
derived from Fan, \etal (1995, p. 147):
\begin{eqnarray} \label{c11}
\hAMPEC(q)= C_{p}(K) \bigg[\frac{a(\psi)\int b''(\theta(x))
q''(m(x)) dx} {\int \{\theta^{(p+1)}(x)\}^2 \{b''(\theta(x))\}^2
q''(m(x)) f_X(x) dx}\bigg]^{1/(2p+3)} n^{-1/(2p+3)},
\end{eqnarray}
where $C_p(K)$ is a constant depending only on the degree and kernel
of the local regression. In particular, for the $q_2$-function which
gives the deviance-loss, we have
\begin{eqnarray} \label{c12}
\hAMPEC(q_2)=C_p(K) \bigg[\frac {a(\psi)|\Omega_X|} {\int
\{\theta^{(p+1)}(x)\}^2 b''(\theta(x)) f_X(x) dx}\bigg]^{1/(2p+3)}
n^{-1/(2p+3)},
\end{eqnarray}
where $|\Omega_X|$ measures the length of the support of $f_X$.
Apparently, this asymptotically optimal bandwidth differs from the
asymptotically optimal bandwidth,
\begin{eqnarray} \label{c13}
\hAMISE =C_p(K) \bigg[\frac {a(\psi)\int \{b''(\theta(x))\}^{-1}
dx} {\int \{\theta^{(p+1)}(x)\}^2 f_X(x) dx}\bigg]^{1/(2p+3)}
n^{-1/(2p+3)},
\end{eqnarray}
determined by minimizing the asymptotic $\MISE(h)$ of $\hat \theta$
defined in (\ref{b9}), with an exception of the Gaussian family.
\begin{center}
{[ \textsl{Put Table \ref{table-1} about here} ]}
\end{center}

The second issue concerns how far away $\hAMPEC(q_2)$ departs from
$\hAMISE$.
 For Poisson and Bernoulli response
variables, examples in Table \ref{table-1} illustrate that the
distinction between $\hAMPEC(q_2)$ and $\hAMISE$ can be
noticeable. To gain further insights, we will need the following
definition.
\begin{definition}
Two functions $F$ and $G$ are called \textsl{similarly ordered} if
$\{F(x_1)-F(x_2)\}\{G(x_1)-G(x_2)\}\ge 0$ for all $x_1$ in the
domain of $F$ and all $x_2$ in the domain of $G$, and
\textsl{oppositely ordered} if the inequality is reversed.
\end{definition}

The following theorem characterizes the relation between
$\hAMPEC(q_2)$ and $\hAMISE$.
\begin{proposition} \label{prop3}
Define $F(x)=\{\theta^{(p+1)}(x)\}^2 b''(\theta(x)) f_X(x)$ and
$G(x)=\{b''(\theta(x))\}^{-1}$. Assume that $p$ is odd.
\begin{itemize}
\item[$(a)$] If $F$ and $G$ are oppositely ordered, then $\hAMPEC(q_2)\le \hAMISE$.
If $F$ and $G$ are similarly ordered, then $\hAMPEC(q_2)\ge
\hAMISE$.

\item[$(b)$] Assume that $b''(\theta(x))$ is bounded away from zero and infinity.
Write $m_{b''}=\min_{x\in \Omega_X} b''(\theta(x))$ and
$M_{b''}=\max_{x\in \Omega_X} b''(\theta(x))$. If $\theta(x)$ is
a polynomial function of degree $p+1$, and $f_X$ is a uniform
density on $\Omega_X$, then
$$
\bigg\{\frac{4m_{b''}
M_{b''}}{(m_{b''}+M_{b''})^2}\bigg\}^{1/(2p+3)} \le
\frac{\hAMPEC(q_2)} {\hAMISE} \le 1,
$$
in which the equalities are satisfied if and only if the exponential
family is Gaussian.
\end{itemize}
\end{proposition}

\subsection{Empirical Cross-Validation} \label{sect-3.4}

The approximate CV criterion (\ref{c9}) can be further simplified.
To this end,  we first approximate the ``degrees of freedom"
$\sum_{i=1}^n H_i$ (Hastie and Tibshirani 1990). To facilitate
presentations, we now define the ``equivalent kernel" $ \calK(t)$
induced by the local-polynomial fitting as the first element of the
vector $S^{-1} (1, t, \ldots, t^{p})^T K(t)$, in which the matrix
$S=(\mu_{i+j-2})_{1\le i, j\le p+1}$ with $\mu_k=\int t^k K(t)\,
dt$. See Ruppert and Wand (1994).
\begin{proposition} \label{prop4}
Assume conditions ${\bf (A)}$ and ${\bf (B)}$ in the Appendix.
If $n \to \infty$,
$h\to 0$, and $nh\to \infty$, we have
\begin{eqnarray*}
\sum_{i=1}^n H_i = \calK(0) |\Omega_X|/h\, \{1+o_P(1)\},
\end{eqnarray*}
where $|\Omega_X|$ denotes the length of the support of the
random variable $X$.
\end{proposition}

Proposition \ref{prop4} shows that the degrees of freedom is
asymptotically independent of the design density and the conditional
density. It approximates the notion of model complexity in
nonparametric fitting.
\begin{center}
{[ \textsl{Put Table \ref{table-2} about here} ]}
\end{center}

Proposition~\ref{prop4} does not specify the constant term. To use
the asymptotic formula for finite samples, we need some bias
corrections. Note that when $h \to \infty$, the local polynomial
fitting becomes a global polynomial fitting. Hence, its degrees of
freedom should be $p+1$. This leads us to propose the following
empirical formula:
\begin{eqnarray} \label{c14}
 \sum_{i=1}^n H_i \doteq (p+1-\ta)+\calC n/(n-1)\calK(0)|\Omega_X|/h.
\end{eqnarray}
In the Gaussian family, Zhang (2003) used simulations to determine
the choices $\ta$ and $\calC$.  See Table \ref{table-2}, which uses
the Epanechnikov kernel function, $K(t)=.75(1-t^2)_+$.
Interestingly, our simulation studies in Section \ref{sect-7}
demonstrate that these choices also work well for Poisson responses.
However, for Bernoulli responses, we find that for $p=1$, slightly
different choices given by $\ta=.7$ and $\calC=1.09$ provide better
approximations.

We propose the empirical version of the estimated total prediction
error by replacing $H_i$ in (\ref{c9}) with their empirical average,
$\overline{H}_E = (p+1-\ta)/n+\calC /(n-1)\calK(0)|\Omega_X|/h$,
leading to the empirical cross-validation (ECV) criterion,
\begin{eqnarray} \label{c15}
\hat \Err^{\ECV}(h) = \sum_{i=1}^n \Big[Q(Y_i,\hat m_{i})
+2^{-1}q''(\hat m_i)(Y_i-\hat m_i)^2 \big\{1-{1}/{(1-\overline{H}_E
)^{2}}\big\} \Big].
\end{eqnarray}
This avoids calculating the smoother matrix $\calH$.  Yet, it turns
out to work reasonably well in practice. A data-driven optimal
bandwidth selector, $\hECV$, can be obtained by minimizing
(\ref{c15}).

\section{Nonparametric Logistic Regression} \label{sect-4}

Nonparametric logistic regression plays a prominent role in
classification and regression analysis. Yet, distinctive challenges
arise from the local MLE and bandwidth selection. When the responses
in a local neighborhood are entirely zeros or entirely ones (or
nearly so), the local MLE does not exist.  M\"uller and Schmitt
(1988, p. 751) reported that the local-likelihood method suffers
from a substantial proportion of ``incalculable estimates".  Fan and
Chen (1999) proposed to add a ridge parameter to attenuate the
problem. The numerical instability problem still exists as the ridge
parameter can be very close to zero.  A numerically viable solution
is the lower bound method, which we now introduce.

\subsection{Lower Bound Method for Local MLE}

The lower-bound method is very simple.  For optimizing a concave
function $\calL$, it replaces the Hessian matrix $\triangledown^2
\calL(\bbeta)$ in the Newton-Raphson algorithm by a negative
definite matrix $\bB$, such that
\begin{eqnarray*}
\triangledown^2 \calL(\bbeta) \ge \bB,\quad \mbox{for all}\
\bbeta.
\end{eqnarray*}
Lemma \ref{lem-1}, shown in B\"ohning (1999, p. 14), indicates that
the Newton-Raphson estimate, with the Hessian matrix replaced by the
surrogate $\bB$, can always enhance the target function $\calL$.
\begin{lemma} \label{lem-1}
Starting from any $\bbeta_0$, the LB iterative estimate, defined by
$\bbeta_{\LB}=\bbeta_0-\bB^{-1}\triangledown \calL(\bbeta_0)$, leads
to a monotonic increase of $\calL(\cdot)$, that is,
$\calL(\bbeta_{\LB})-\calL(\bbeta_0)\ge -2^{-1} \triangledown
\calL(\bbeta_0)^T \bB^{-1} \triangledown \calL(\bbeta_0) \ge 0$.
\end{lemma}

For the local logistic regression, $\triangledown^2 \ell(\bbeta; x)
= - \bX(x)^T \bW(x;\bbeta) \bX(x)$. Since $\bzero \le \bW(x;\bbeta)
\le 4^{-1}\bK(x)$, where $\bK(x)=\diag\{K_h(X_j-x)\}$, the Hessian
matrix $\triangledown^2 \ell(\bbeta; x)$ indeed has a lower bound,
\begin{eqnarray} \label{d1}
\bB(x)=-4^{-1} \bX(x)^T \bK(x) \bX(x),
\end{eqnarray}
and the LB-adjusted Newton-Raphson iteration for computing $\hat
\bbeta(x)$ becomes
\begin{eqnarray} \label{d2}
{\bbeta}_{L} = {\bbeta}_{L-1}  - \{\bB(x)\}^{-1} \bX(x)^T \bK(x)
\br(x; {\bbeta}_{L-1}), \quad L=1,2,\ldots,
\end{eqnarray}
where $\br(x;\bbeta)=(Y_1-m_1(x;\bbeta),\ldots,Y_n-m_n(x;\bbeta))^T$
with $m_j(x;\bbeta)={1}/{[1+\exp\{-\bx_j(x)^T \bbeta\}]}$.

The LB method offers a number of advantages to compute $\hat
\bbeta(x)$. Firstly, the corresponding LB matrix $\bB(x)$ is free
of the parameter vector $\bbeta$, and thus can be computed in
advance of the NR iteration. This in turn reduces the
computational cost. Secondly, the LB matrix is stable, as it is
the same matrix used in the least-squares local-polynomial
regression estimates and does not depend on estimated local
parameters. Thirdly, since the local-likelihood function
$\ell(\bbeta; x)$ is concave, the LB iteration is guaranteed to
increase $\ell(\bbeta; x)$ at each step and converge to its global
maximum $\hat \bbeta(x)$.

\subsection{A Hybrid Bandwidth Selection Method} \label{sect-4.2}

For binary responses, our simulation studies show that the bandwidth
choice minimizing (\ref{c9}) or its empirical version (\ref{c15})
tends to produce over-smoothed estimates. Such a problem was also
encountered in Aragaki and Altman (1997) and Fan, Farmen and Gijbels
(1998, Table 1). Because of the importance of binary responses in
nonparametric regression and classification,  a new bandwidth
selector that specifically accommodates the binary responses is
needed.

We first employ the LB scheme (\ref{d2}) to derive a new one-step
estimate of $\hat \bbeta^{-i}(x)$, starting from $\hat \bbeta(x)$.
Define $S_n(x)=\bX(x)^T \bK(x)\bX(x)$ and
$\calS_i=\bevector_1^T\{S_n(X_i)\}^{-1}\bevector_1K_h(0)$, where
$\bevector_1=(1,0,\ldots,0)^T$. The resulting leave-one-out formulas
and the cross-validated estimates of the total prediction error are
displayed in Proposition \ref{prop5}.
\begin{proposition} \label{prop5}
Assume conditions $(A1)$ and $(A2)$ in the Appendix. Then for the
local-likelihood MLE in the Bernoulli family,
\begin{eqnarray}
\hat \bbeta_{LB}^{-i}(x)-\hat \bbeta(x) &\doteq&
-\frac{4\{S_n(x)\}^{-1}\bx_i(x)K_h(X_i-x)\{Y_i-b'(\bx_i(x)^T\hat
\bbeta(x))\}}
         {1-K_h(X_i-x)\bx_i(x)^T\{S_n(x)\}^{-1}\bx_i(x)}, \label{d3}\\
\hat \theta_i^{-i}-\hat \theta_i
&\doteq& -\frac{4\calS_i}{1-\calS_i} (Y_i-\hat m_i), \label{d4}\\
\hat m_i^{-i}-\hat m_i
&\doteq& -\frac{4 b''(\hat \theta_i)\calS_i}{1-\calS_i} (Y_i-\hat m_i), \label{d5} \\
{\hat \Err}^{\CV} &\doteq& \sum_{i=1}^n \Big[Q(Y_i,\hat
m_{i})+2^{-1} q''(\hat m_i) (Y_i-\hat m_i)^2
         \big[1-\{1+{4b''(\hat \theta_i)\calS_i}/{(1-\calS_i)}\}^2\big]\Big].   \label{d6}
\end{eqnarray}
\end{proposition}

Direct use of a bandwidth selector that minimizes (\ref{d6}) tends
to under-smooth the binary responses. To better appreciate this,
note that the second term in (\ref{d6}) is approximately
\begin{eqnarray} \label{d7}
 -q''(\hat m_i) (Y_i-\hat m_i)^2 \{4b''(\hat
\theta_i)\}\calS_i,
\end{eqnarray}
and the second in (\ref{c9}) can be approximated as
\begin{eqnarray} \label{d8}
-q''(\hat m_i) (Y_i-\hat m_i)^2 H_i.
\end{eqnarray}
As demonstrated in Lemma \ref{lem-3} in the Appendix, $\calS_i$
decreases with $h$  and $H_i\doteq \calS_i$. Since $0\le 4b''(\hat
\theta_i) \le 1$ for the Bernoulli family, (\ref{d7}) down weighs
the effects of model complexity, resulting in a smaller bandwidth.

The above discussion leads us to define a hybrid version of ${\hat
\Err}^{\CV}$ as
\begin{eqnarray}
\sum_{i=1}^n \Big[Q(Y_i,\hat m_{i})+2^{-1} q''(\hat m_i) (Y_i-\hat
m_i)^2 \big[1-\{1+{2b''(\hat
\theta_i)\calS_i}/{(1-\calS_i)}+{2^{-1}H_i}/{(1-H_i)}\}^2
\big]\Big],
  \label{d9}
\end{eqnarray}
which averages terms in (\ref{d7}) and (\ref{d8}) to mitigate the
oversmoothing problem of criterion (\ref{c9}). This new criterion
has some desirable properties: ${2b''(\hat
\theta_i)\calS_i}/{(1-\calS_i)}+{2^{-1}H_i}/{(1-H_i)}$ is bounded
below by ${2^{-1}H_i}/{(1-H_i)}$, thus guarding against
under-smoothing, and is bounded above by
$\{{\calS_i}/{(1-\calS_i)}+{H_i}/{(1-H_i)}\}/2$, thus diminishing
the influence of over-smoothing. An empirical cross-validation
criterion is to replace $\calS_i$ and $H_i$ in (\ref{d9}) by their
empirical averages, which are (\ref{c14}) divided by $n$. A hybrid
bandwidth selector for binary responses can be obtained by
minimizing this ECV.

\section{Extension to Generalized Varying-Coefficient Model} \label{sect-5}

This section extends the techniques of Sections \ref{sect-3} and
\ref{sect-4} to a useful class of multi-predictor models. The
major results are presented in Propositions
\ref{prop6}--\ref{prop8}.

Consider multivariate predictor variables, containing a scalar $U$
and a vector $\bvectorX=(X_1, \ldots, X_d)^T$. For the response
variable $Y$ having a distribution in the exponential-family, define
by $m(u, \bvectorx)=E(Y|U=u, \bvectorX=\bvectorx)$ the conditional
mean regression function, where $\bvectorx=(x_1, \ldots, x_d)^T$.
The generalized varying-coefficient model assumes that the
$d+1$-variate canonical parameter function
$\theta(u,\bvectorx)=g(m(u,\bvectorx))$, with the canonical link
$g$, takes the form
\begin{eqnarray} \label{e1}
g(m(u, \bvectorx)) = \theta(u,\bvectorx)=\sum_{k=1}^{d}
a_k(u)x_k=\bvectorx^T \bA(u).
\end{eqnarray}
for a vector $\bA(u)=(a_1(u) , \ldots, a_d(u))^T$ of unknown smooth
coefficient functions.

We first describe the local-likelihood estimation of $\bA(u)$, based
on the independent observations $\{(U_j, \bvectorX_j,
Y_j)_{j=1}^{n}\}$. Assume that $a_k(\cdot)$'s are $(p+1)$-times
continuously differentiable at a fitting point $u$. Put
$\bA^{(\ell)}(u)=(a_1^{(\ell)}(u) , \ldots, a_d^{(\ell)}(u))^T$.
Denote by $\bbeta(u)=(\bA(u)^T, \bA^{(1)}(u)^T, \ldots,
\bA^{(p)}(u)^T/p!)^T$ the $d(p+1)$ by $1$ vector of coefficient
functions along with their derivatives, $\bu_j(u)=(1, (U_j-u),
\ldots, (U_j-u)^p)^T$, and $\bI_d$ a $d\times d$ identity matrix.
For observed covariates $U_j$ close to the point $u$,
\begin{eqnarray*}
\bA(U_j) \doteq
\bA(u)+(U_j-u)\bA^{(1)}(u)+\cdots+(U_j-u)^p\bA^{(p)}(u)/p!
=\{\bu_j(u) \otimes \bI_d\}^T \bbeta(u),
\end{eqnarray*}
in which the symbol $\otimes$ denotes the Kronecker product, and
thus from (\ref{e1}),
\begin{eqnarray*}
\theta(U_j,\bvectorX_j) \doteq \{\bu_j(u)\otimes \bvectorX_j\}^T
\bbeta(u).
\end{eqnarray*}
The local-likelihood MLE $\hat \bbeta(u)$ maximizes the local
log-likelihood function:
\begin{eqnarray} \label{e2}
\ell(\bbeta; u) = \sum_{j=1}^n l(\{\bu_j(u)\otimes \bvectorX_j\}^T
\bbeta; Y_j) K_h(U_j-u).
\end{eqnarray}
The first $d$ entries of $\hat \bbeta(u)$ supply the local MLEs
$\hat \bA(u)$ of $\bA(u)$, and the local MLEs of
$\theta(u,\bvectorx)$ and $m(u,\bvectorx)$ are given by $\hat
\theta(u,\bvectorx)=\bvectorx^T \hat \bA(u)$ and $\hat
m(u,\bvectorx)=b'(\hat \theta(u,\bvectorx))$, respectively. A
similar estimation procedure, applied to $n-1$ observations
excluding $(U_i, \bvectorX_i, Y_i)$, leads to the local
log-likelihood function, $\ell^{-i}(\bbeta; u)$, and the
corresponding local MLEs, $\hat \bbeta^{-i}(u)$, $\hat
\theta^{-i}(u,\bvectorx)$, and $\hat m^{-i}(u,\bvectorx)$
respectively.

\subsection{Leave-One-Out Formulas} \label{sect-5.1}

To derive the leave-one-out formulas in the case of multivariate
covariates, we need some additional notations. Let
$\bX^*(u)=(\bu_1(u) \otimes \bvectorX_1,\ldots,\bu_n(u) \otimes
\bvectorX_n)^T$, $\bW^*(u;\bbeta)=\diag\{K_h(U_j-u) b''(\{\bu_j(u)
\otimes \bvectorX_j\}^T\bbeta)\}$, and $S_n^*(u;\bbeta)=\bX^*(u)^T
\bW^*(u;\bbeta) \bX^*(u)$. Define a projection matrix as
\begin{eqnarray*}
\calH^*(u;\bbeta)=\{\bW^*(u;\bbeta)\}^{1/2} \bX^*(u)
\{S_n^*(u;\bbeta)\}^{-1} \bX^*(u)^T \{\bW^*(u;\bbeta)\}^{1/2}.
\end{eqnarray*}
Let $\calH^*_{ii}(u;\bbeta)$ be its $i$th diagonal entry and
$H_i^*=\calH^*_{ii}(U_i;\hat \bbeta(U_i))$. Propositions \ref{prop6}
and \ref{prop7} below present the leave-one-out formulas and
cross-validated estimate of the total prediction error.
\begin{proposition} \label{prop6}
Assume condition $(A2)$ in the Appendix. Then for $i=1, \ldots, n$,
\begin{eqnarray*}
\hat \bbeta^{-i}(u)-\hat \bbeta(u) &\doteq& -\frac
{\{S_n^*(u;\hat \bbeta(u))\}^{-1} \{\bu_i(u) \otimes \bvectorX_i\}
         K_h(U_i-u) \{Y_i-b'(\{\bu_i(u) \otimes \bvectorX_i\}^T \hat \bbeta(u)) \}}
         {1-\calH^*_{ii}(u;\hat \bbeta(u))},  \\
\hat \theta^{-i}_i-\hat \theta_i
&\doteq& -\frac {H_i^*} {1-H_i^*} (Y_i-\hat m_i)/b''(\hat \theta_i),  \\
\hat m^{-i}_i-\hat m_i &\doteq& -\frac {H_i^*} {1-H_i^*}
(Y_i-\hat m_i).
\end{eqnarray*}
\end{proposition}
\begin{proposition} \label{prop7}
Assume conditions $(A1)$ and $(A2)$ in the Appendix. Then
\begin{eqnarray}
{\hat \Err}^{\CV} \doteq \sum_{i=1}^n \Big[Q(Y_i,\hat m_{i}) +2^{-1}
q''(\hat m_i)(Y_i-\hat m_i)^2 \big\{1-{1}/{(1-H_i^*)^{2}}\big\}
\Big]. \label{e3}
\end{eqnarray}
\end{proposition}

\subsection{Empirical Cross-Validation} \label{sect-5.2}

In the generalized varying-coefficient model, the asymptotic
expression of the degrees of freedom $\sum_{i=1}^n H_i^*$ is given
below.
\begin{proposition} \label{prop8}
Assume conditions $({\bf A})$ and $({\bf C})$ in the Appendix.
If $n \to \infty$,
$h\to 0$, and $nh\to \infty$, we have
\begin{eqnarray*}
\sum_{i=1}^n H_i^* = d\calK(0) |\Omega_U|/h\, \{1+o_P(1)\}.
\end{eqnarray*}
\end{proposition}

As $h\to \infty$, the total number of model parameters becomes
$d(p+1)$ and this motivates us to propose the empirical formula
for degrees of freedom:
\begin{eqnarray} \label{e4}
\sum_{i=1}^n H_i^* \doteq d \{(p+1-\ta)+\calC n/(n-d)\calK(0)
|\Omega_U|/h\}.
\end{eqnarray}
The empirical version of the estimated total prediction error is to
replace $H_i^*$ in (\ref{e3}) by $d \{(p+1-\ta)/n+\calC
/(n-d)\calK(0) |\Omega_U|/h\}$. Call ${\hat \Err}^{\ECV}(h)$ the
empirical version of the cross-validation criterion. Compared with
the bandwidth selector in Cai, Fan and Li (2000), the ${\hat
\Err}^{\ECV}(h)$-minimizing bandwidth selector, $\hECV$, is much
easier to obtain.

\subsection{Binary Responses} \label{sect-5.3}

For Bernoulli responses, the LB method in Section \ref{sect-4}
continues to be applicable for obtaining $\hat \bbeta(u)$ and $\hat
\bbeta^{-i}(u)$. For the local logistic regression, $\triangledown^2
\ell(\bbeta; u)$ has a lower bound, $\bB(u)=- 4^{-1}\bX^*(u)^T
\bK^*(u)\bX^*(u)$, where $\bK^*(u)=\diag\{K_h(U_j-u)\}$. Similar to
(\ref{d2}), the LB-adjusted NR iteration for $\hat \bbeta(u)$
proceeds as follows,
\begin{eqnarray*}
\bbeta_{L}=\bbeta_{L-1}-\{\bB(u)\}^{-1} \bX^*(u)^T \bK^*(u)
\br^*(u;\bbeta_{L-1}), \quad L=1,2,\ldots,
\end{eqnarray*}
where
$\br^*(u;\bbeta)=(Y_1-m_1^*(u;\bbeta),\ldots,Y_n-m_n^*(u;\bbeta))^T$
with $m_j^*(u;\bbeta)={1}/{(1+\exp[-\{\bu_j(u) \otimes
\bvectorX_j\}^T \bbeta])}$. The leave-one-out formulas and the
cross-validated estimates of the prediction error are similar to
those in Proposition \ref{prop5}, with $S_n(x)$ replaced by
$S_n^*(u)=\bX^*(u)^T \bK^*(u)\bX^*(u)$ and $\calS_i$ by
$\calS_i^*=(\bevector_1 \otimes \bvectorX_i)^T\{S_n^*(U_i)\}^{-1}
(\bevector_1 \otimes \bvectorX_i)K_h(0)$. In the spirit of
(\ref{d9}), the hybrid selection criterion for bandwidth is
\begin{eqnarray}
\sum_{i=1}^n \Big[Q(Y_i,\hat m_{i})+2^{-1} q''(\hat m_i) (Y_i-\hat
m_i)^2 \big[1-\{1+{2b''(\hat
\theta_i)\calS_i^*}/{(1-\calS_i^*)}+{2^{-1}H_i^*}/{(1-H_i^*)}\}^2
\big]\Big].   \label{e5}
\end{eqnarray}
The ECV criterion can be obtained similarly via replacing
$\calS_i^*$ and $H_i^*$ by their empirical averages, which are
(\ref{e4}) divided by $n$.

\section{Extension to Semiparametric Model} \label{sect-6}

A further extension of model (\ref{e1}) is to allow part of the
covariates independent of $U$, resulting in
\begin{equation}
\theta(u,\bvectorx, \bvectorz)= \bvectorx^T \bA(u) + \bvectorz^T
\bbeta, \label{f1}
\end{equation}
in which $(u,\bvectorx, \bvectorz)$ lists values of all covariates
$(U,\bvectorX, \bvectorZ)$. This model keeps the flexibility that
$\bvectorZ$ does not interact with $U$.  The challenge is how to
choose a bandwidth to efficiently estimate both the parametric and
nonparametric components. In this section, we propose a two-stage
bandwidth selection method which is applicable to general
semiparametric models. This is an attempt to answer an important
question raised by Bickel and Kwon (2001) on the bandwidth selection
for semiparametric models.

The parameters in model (\ref{f1}) can be estimated via the profile
likelihood method. For each given $\bbeta$,  applying the
local-likelihood method with a bandwidth $h$, we obtain an estimate
$\hat{\bA}(u; \bbeta, h)$, depending on $h$ and $\bbeta$.
Substituting it into (\ref{f1}), we obtain a pseudo parametric
model:
\begin{equation}
\theta(u,\bvectorx, \bvectorz) \doteq \bvectorx^T \hat{\bA}(u;
\bbeta, h) + \bvectorz^T \bbeta. \label{f2}
\end{equation}
Regarding (\ref{f2}) as a parametric model with parameter $\bbeta$,
by using the maximum likelihood estimation method, we obtain the
profile likelihood estimators $\hat{\bbeta}(h)$ and $\hat{\bA}(u;
\hat{\bbeta}(h), h)$. This estimator is semiparametrically
efficient.

We now outline a two-stage method for choosing the bandwidth. The
idea is also applicable to other semiparametric problems. Starting
from a very small bandwidth $h_0$, we obtain a semiparametric
estimator $\hat{\bbeta}(h_0)$ (see a justification below).  To avoid
difficulty of implementation, the nearest type of bandwidth can be
used. This estimator is usually root-n consistent. Thus, $\bbeta$ in
model (\ref{f1}) can be regarded as known and hence model (\ref{f1})
becomes a varying coefficient model. Applying a bandwidth selection
method for varying coefficient models, such as the approximate
cross-validation method in the previous section, we obtain a
bandwidth $\hat{h}$. Using this $\hat{h}$, we obtain the profile
likelihood estimator $\hat{\bbeta}(\hat h)$ and $\hat{\bA}(u;
\hat{\bbeta}(\hat h), \hat h)$.   This is a two-stage method for
choosing the bandwidth for a semiparametric model.

To illustrate the idea, we specifically consider the partially
linear model:
\begin{equation}
   Y_i = a(U_i) + \bvectorZ_i^T \bbeta + \varepsilon_i, \quad i= 1,\ldots, n.
   \label{f3}
\end{equation}
Assume that the data have already been sorted according to $U_i$.
Let $h_0$ be the nearest two-point bandwidth so that
$$
  \hat{a}(U_i; \bbeta, h_0) = 2^{-1} (Y_i - \bvectorZ_i^T \bbeta +
  Y_{i-1} - \bvectorZ_{i-1}^T \bbeta).
$$
Substituting this into (\ref{f3}) and rearranging the equation, we
have
\begin{equation}
  Y_i - Y_{i-1} \doteq   (\bvectorZ_i  - \bvectorZ_{i-1})^T \bbeta + 2\varepsilon_i.
  \label{f4}
\end{equation}
Applying the least-squares method, we obtain an estimator
$\hat{\bbeta}(h_0)$.

To see why such a crude parametric estimator $\hat{\bbeta}(h_0)$
is root-n consistent, let us take the difference of (\ref{f3}).
Under the mild conditions, $a(U_i) - a(U_{i-1}) = O_P(n^{-1})$.
Hence, the difference of (\ref{f3}) yields
\[
Y_i - Y_{i-1} =  (\bvectorZ_i  - \bvectorZ_{i-1})^T \bbeta +
\varepsilon_i - \varepsilon_{i-1} + O_P(n^{-1}).
\]
Hence, the least-squares estimator  $\hat{\bbeta}$, which is the
same as $\hat{\bbeta}(h_0)$, is root-n consistent.  See Yachew
(1997) for a proof.  This example shows that even if a very crude
bandwidth $h_0$ is chosen, the parametric component
$\hat{\bbeta}(h_0)$ is still root-n consistent.

The two-stage bandwidth selector is to apply a data-driven
bandwidth selector to the following univariate nonparametric
regression problem:
$$
Y_i = a(U_i) + \bvectorZ_i^T \hat{\bbeta}(h_0) + \varepsilon_i,
$$
and to use the resulting bandwidth for the original semiparametric
problem.  Such an idea was implemented in  Fan and Huang (2005).
They reported that the resulting semiparametric and nonparametric
estimators are efficient.

\section{Simulations} \label{sect-7}

The purpose of the simulations is three-fold: to assess the accuracy
of the empirical formulas (\ref{c14}) and (\ref{e4}), the
performance of the bandwidth selector $\hECV$, and the behavior of
the proposed bandwidth selector for local-likelihood estimation. For
Bernoulli responses, we apply the hybrid bandwidth selector to local
logistic regression. Throughout our simulations, we use the
$q_2$-function associated with the deviance loss for bandwidth
selection, combined with the local-linear likelihood method and the
Epanechnikov kernel. Unless specifically mentioned otherwise, the
sample size is $n=400$. A complete copy of Matlab codes is available
upon request.

\subsection{Generalized Nonparametric Regression Model}  \label{sect-7.1}

\begin{figure}[htbp]
\centerline{\psfig{figure=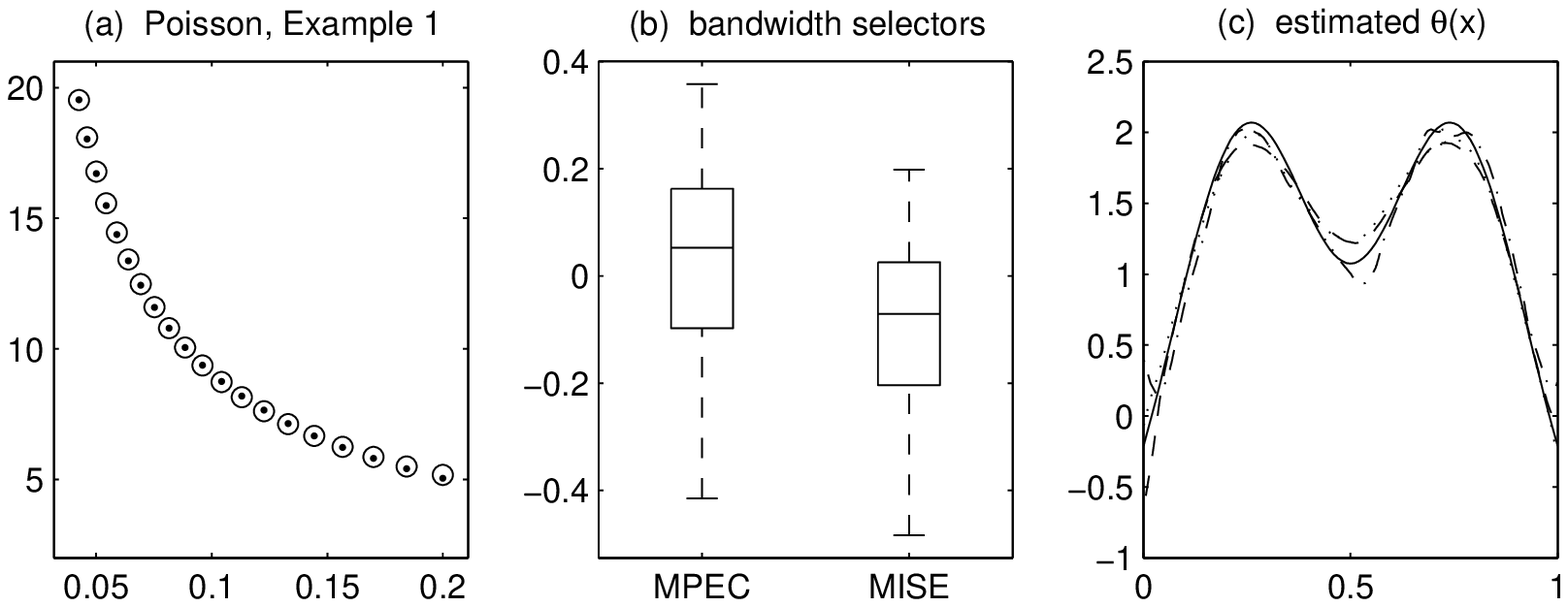,height=2in,width=6in}}
\vskip 0.3in
\centerline{\psfig{figure=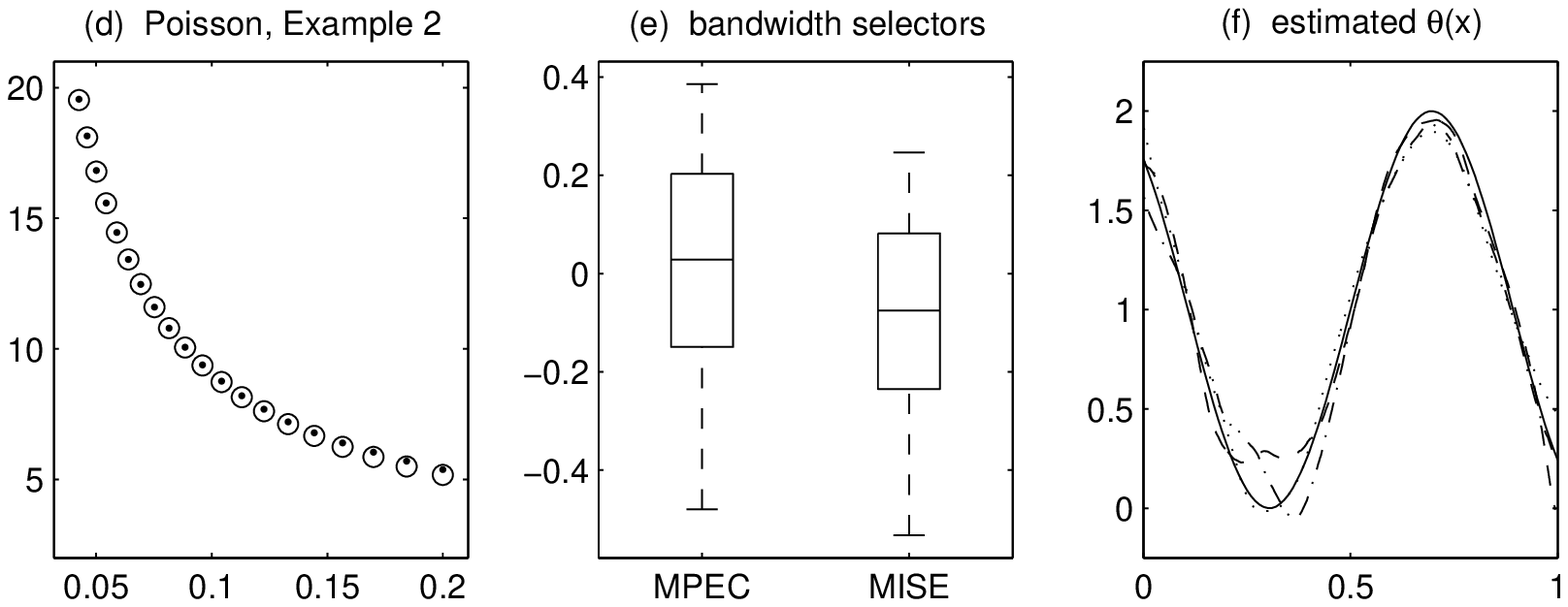,height=2in,width=6in}}
\vskip 0.3in
\centerline{\psfig{figure=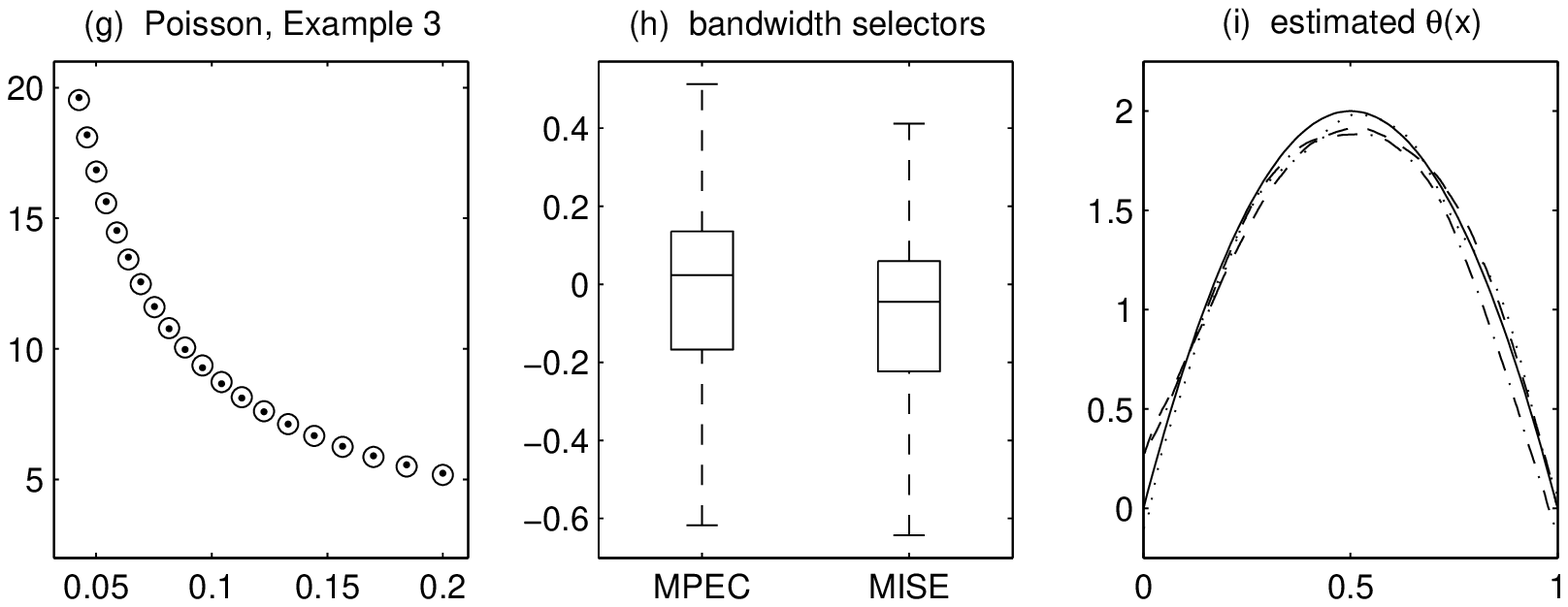,height=2in,width=6in}}
\begin{singlespace}
\caption {\textsl{Evaluation of Local-Likelihood Nonparametric
Regression for Poisson Responses. Left panel: plots of $\sum_{i=1}^n
H_i$ versus $h$. Dots denote the actual values, centers of circles
stand for the empirical values given by (\ref{c14}), for
local-linear smoother with $\ta=.70$ and $\calC=1.03$. Middle panels
[figures (b), (e), and (h)]: boxplots of
$\{\hECV-\hAMPEC(q_2)\}/\hAMPEC(q_2)$ and
$\{\hECV-\hAMISE\}/\hAMISE$. Panels (c), (f), and (i): estimated
curves from three typical samples are presented corresponding to the
$25^{th}$ (the dotted curve), the $50^{th}$ (the dashed curve), and
the $75^{th}$ (the dash-dotted curve) percentiles among the
ASE-ranked values. The solid curves denote the true functions.}}
\label{figure-3}
\end{singlespace}
\end{figure}

For simplicity, we assume that the predictor variable $X$ has the
uniform probability density on the interval $(0,1)$. The bandwidth
$\hECV$ is searched over an interval, $[\hmin, .5]$, at a geometric
grid of $30$ points. We take $\hmin=3\, h_0$ for Poisson regression,
whereas for logistic regression, we take $\hmin=5\, h_0$ in Example
1 and $\hmin=.1$ in Examples 2--3, where $h_0=\max[5/n, \max_{2\le
j\le n} \{X_{(j)}-X_{(j-1)}\}]$, with $X_{(1)} \le \cdots \le
X_{(n)}$ being the order statistics.
\begin{center}
{[ \textsl{Put Figure \ref{figure-3} about here} ]}
\end{center}

\begin{figure}[htbp]
\centerline{\psfig{figure=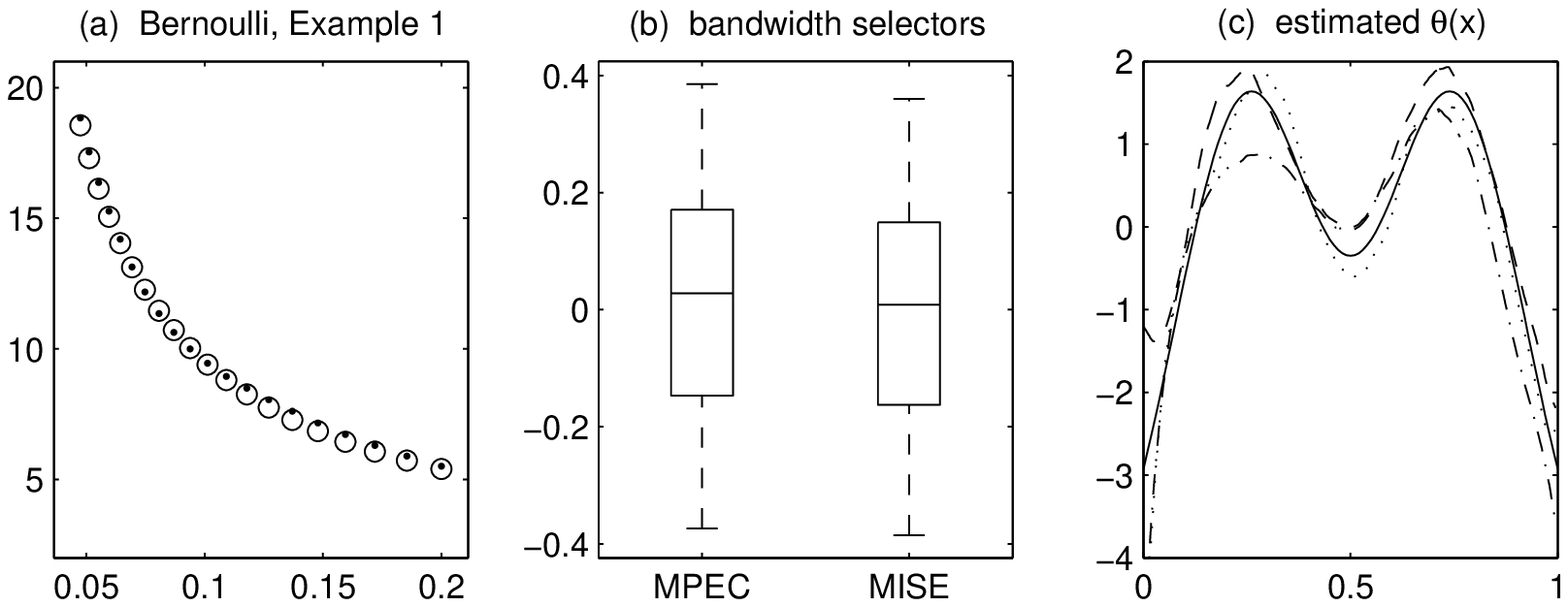,height=2in,width=6in}}
\vskip 0.3in
\centerline{\psfig{figure=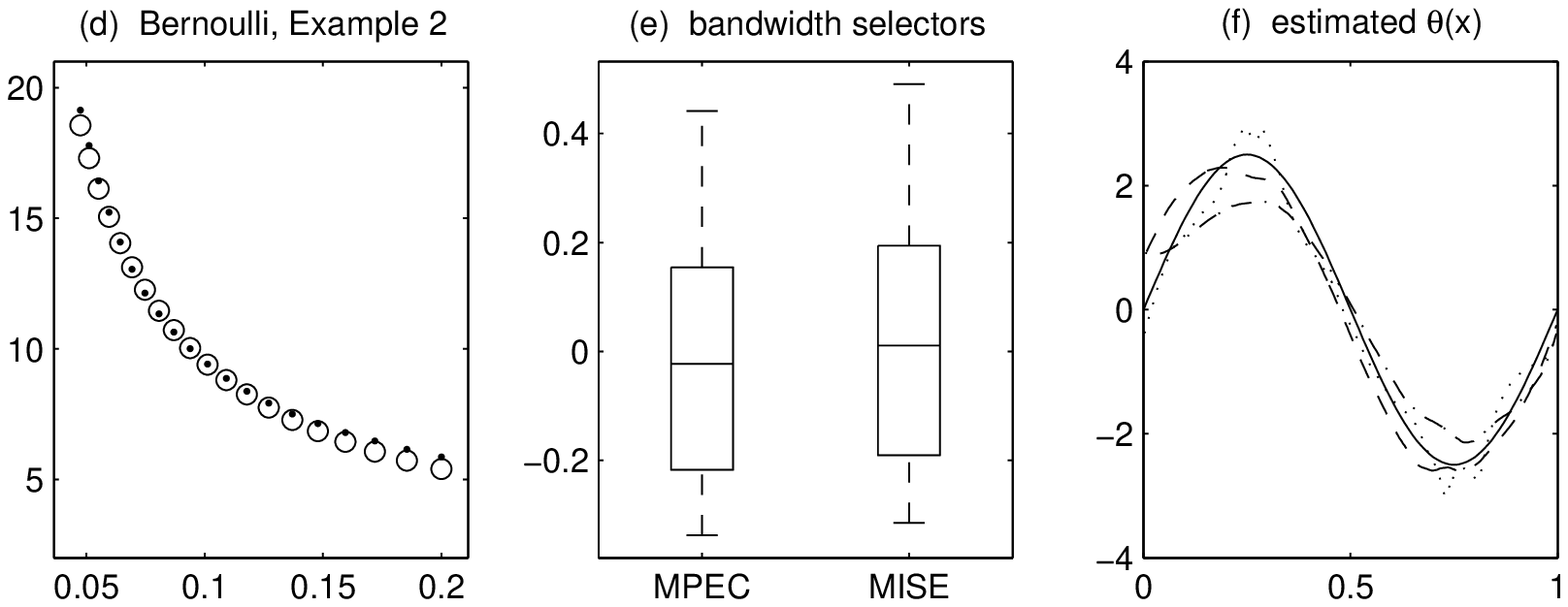,height=2in,width=6in}}
\vskip 0.3in
\centerline{\psfig{figure=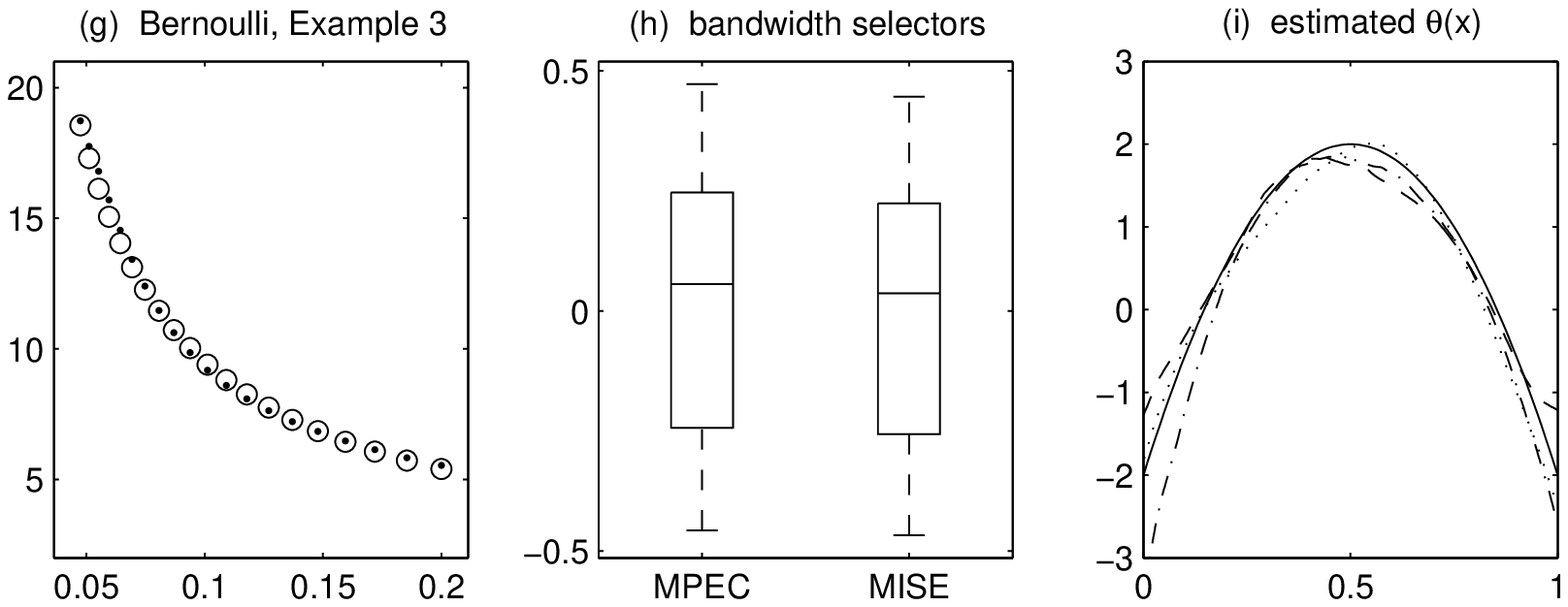,height=2in,width=6in}}
\begin{singlespace}
\caption {\textsl{Evaluation of Local-Likelihood Nonparametric
Regression for Bernoulli Responses. Captions are similar to those
for Figure \ref{figure-3}. Here $\hECV$ minimizes the empirical
version of (\ref{d9}); the empirical formula (\ref{c14}) uses
$\ta=.70$ and $\calC=1.09$ for $H_i$ and $\ta=.70$ and $\calC=1.03$
for $\calS_i$.}} \label{figure-4}
\end{singlespace}
\end{figure}

{\bf Poisson regression:} We first consider the response variable
$Y$ which, conditional on $X=x$, follows a Poisson distribution
with parameter $\lambda(x)$:
$$
P(Y=y|X=x)=\{\lambda(x)\}^y \exp\{-\lambda(x)\}/y!,\quad
y=0,1,\ldots.
$$
The function $\theta(x)=\ln\{\lambda(x)\}$ is given in the test
examples,
\begin{eqnarray*}
\mbox{Example 1:} && \theta(x)=3.5[\exp\{-(4x-1)^2\}+\exp\{-(4x-3)^2\}]-1.5, \\
\mbox{Example 2:} && \theta(x)=\sin\{2(4x-2)\}+1.0, \\
\mbox{Example 3:} && \theta(x)=2-.5(4x-2)^2.
\end{eqnarray*}
As an illustration, we first generate from $(X,Y)$ one sample of
independent observations $\{(X_j, Y_j)_{j=1}^n\}$. Figure
\ref{figure-3}(a) plots the degrees of freedom as a function of $h$.
It is clearly seen that the actual values (denoted by dots) are well
approximated by the empirical values (denoted by circles) given by
(\ref{c14}). To see the performance of $\hECV$, Figure
\ref{figure-3}(b) gives boxplots of the relative error,
$\{\hECV-\hAMPEC(q_2)\}/\hAMPEC(q_2)$ and
$\{\hECV-\hAMISE\}/\hAMISE$, based on $100$ random samples; refer to
Table \ref{table-1} for values of $\hAMPEC(q_2)$ and $\hAMISE$. We
observe that $\hECV$ is closer to $\hAMPEC(q_2)$ than to $\hAMISE$;
this is in accordance with the discussion of Section \ref{sect-3.3}.
In  Figure \ref{figure-3}(c), we simulate another $100$ random
samples and for each set obtain $\hECV$ to estimate $\theta(x)$. We
present the estimated curves from three typical samples. The typical
samples are selected in such a way that their ASE values, in which
$\mbox{ASE}=n^{-1}\sum_{j=1}^{n} \{\hat
\theta(X_j)-\theta(X_j)\}^2$, are equal to the $25$th (dotted line),
$50$th (dashed line), and $75$th (dash-dotted line) percentiles in
the $100$ replications. Inspection of these fitted curves suggests
that the bandwidth selector based on minimizing the cross-validated
deviance does not exhibit undersmoothing in the local-likelihood
regression estimation. In Figure \ref{figure-3}, similar results are
also displayed in the middle panel [Figures \ref{figure-3}(d)--(f)]
for Example 2, and in the bottom panel [Figures
\ref{figure-3}(g)--(i)] for Example 3.
\begin{center}
{[ \textsl{Put Figure \ref{figure-4} about here} ]}
\end{center}

{\bf Logistic regression:} We now consider the Bernoulli response
variable $Y$ with canonical parameter,
$\theta(x)=\logit\{P(Y=1|X=x)\}$, chosen according to
\begin{eqnarray*}
\mbox{Example 1:} && \theta(x)=7[\exp\{-(4x-1)^2\}+\exp\{-(4x-3)^2\}]-5.5, \\
\mbox{Example 2:} && \theta(x)=2.5\, \sin(2\pi x), \\
\mbox{Example 3:} && \theta(x)=2-(4x-2)^2.
\end{eqnarray*}
In Figure \ref{figure-4}, we conduct the simulation experiments
serving a similar purpose to Figure \ref{figure-3}. Plots in the
middle (vertical) panel lend support to the convergence of the
hybrid bandwidth selector $\hECV$ to $\hAMPEC(q_2)$, without
suffering from the under- or over-smoothing problem.

\subsection{Generalized Varying-Coefficient Model} \label{sect-7.2}

\begin{figure}[htbp]
\centerline{\psfig{figure=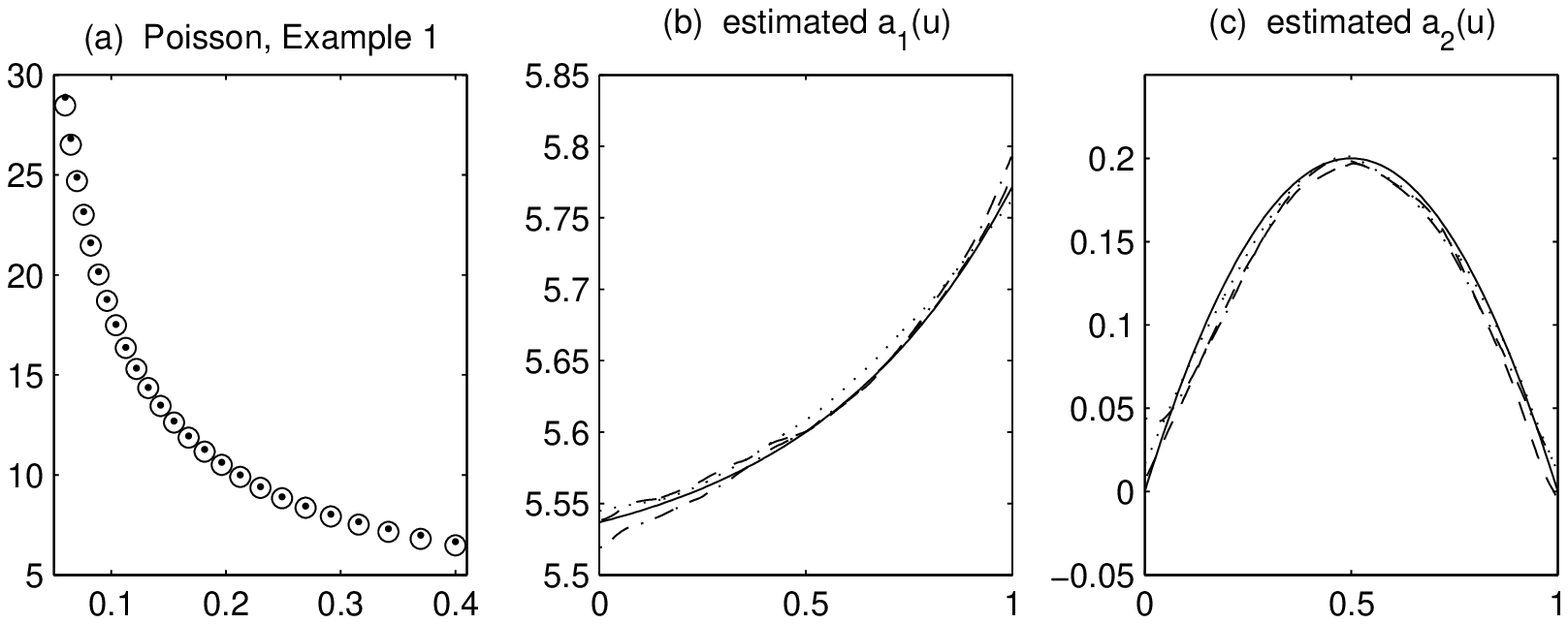,height=2in,width=6in}}
\vskip 0.3in
\centerline{\psfig{figure=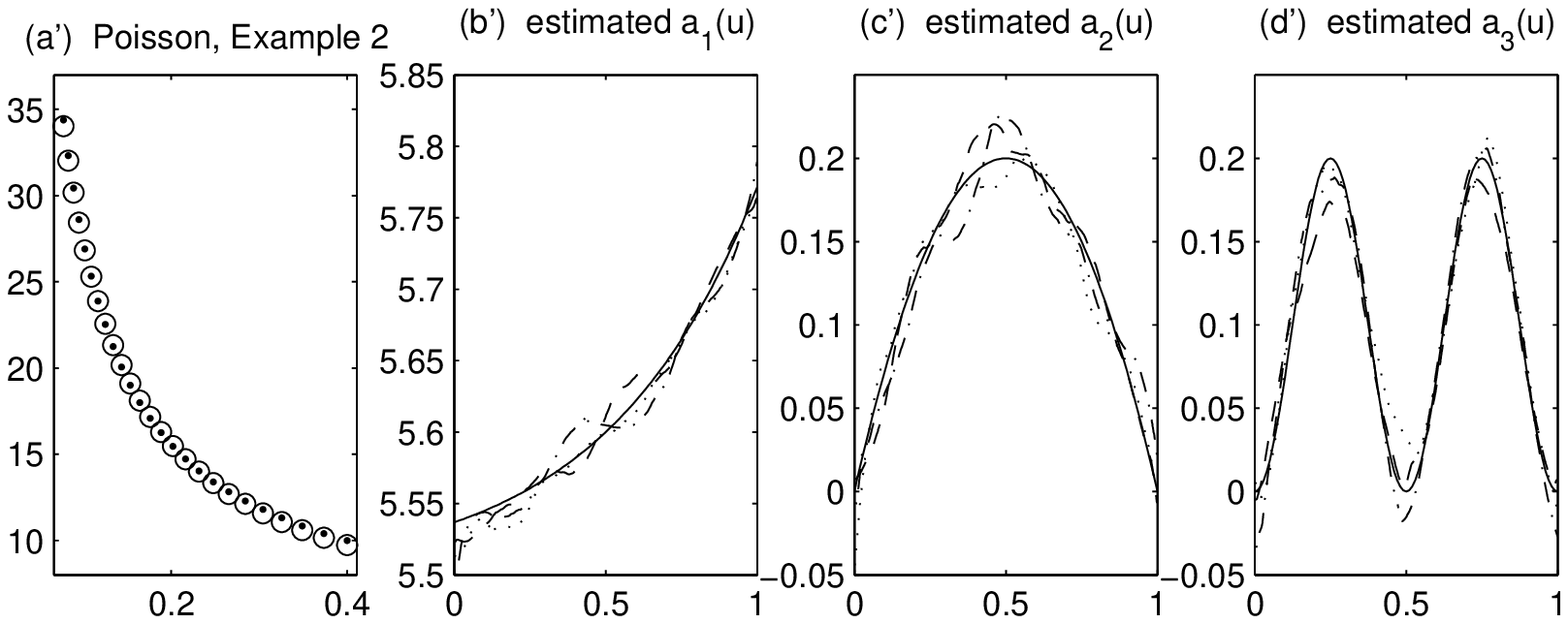,height=2in,width=6in}}
\begin{singlespace}
\caption {\textsl{Evaluation of Local-Likelihood Varying Coefficient
Regression for Poisson Responses. Panels (a) and (a'): plots of
$\sum_{i=1}^n H_i^*$ versus $h$. Dots denote the actual values,
centers of circles stand for the empirical values given by
(\ref{e4}), for local-linear smoother with $\ta=.70$ and
$\calC=1.03$. Panels (b)-(c) for Example 1 and  panels (b')-(d') for
Example 2: estimated curves from three typical samples are presented
corresponding to the $25^{th}$ (the dotted curve), the $50^{th}$
(the dashed curve), and the $75^{th}$ (the dash-dotted curve)
percentiles among the ASE-ranked values. The solid curves denote the
true functions.}} \label{figure-5}
\end{singlespace}
\end{figure}

We consider examples of the generalized varying-coefficient model
(\ref{e1}). We take $\hmin=3\, h_0$ for Poisson regression, where
$h_0=\max[5/n, \max_{2\le j\le n} \{U_{(j)}-U_{(j-1)}\}]$, and
$\hmin=.1$ for logistic regression.
\begin{center}
{[ \textsl{Put Figure \ref{figure-5} about here} ]}
\end{center}

\begin{figure}[htbp]
\centerline{\psfig{figure=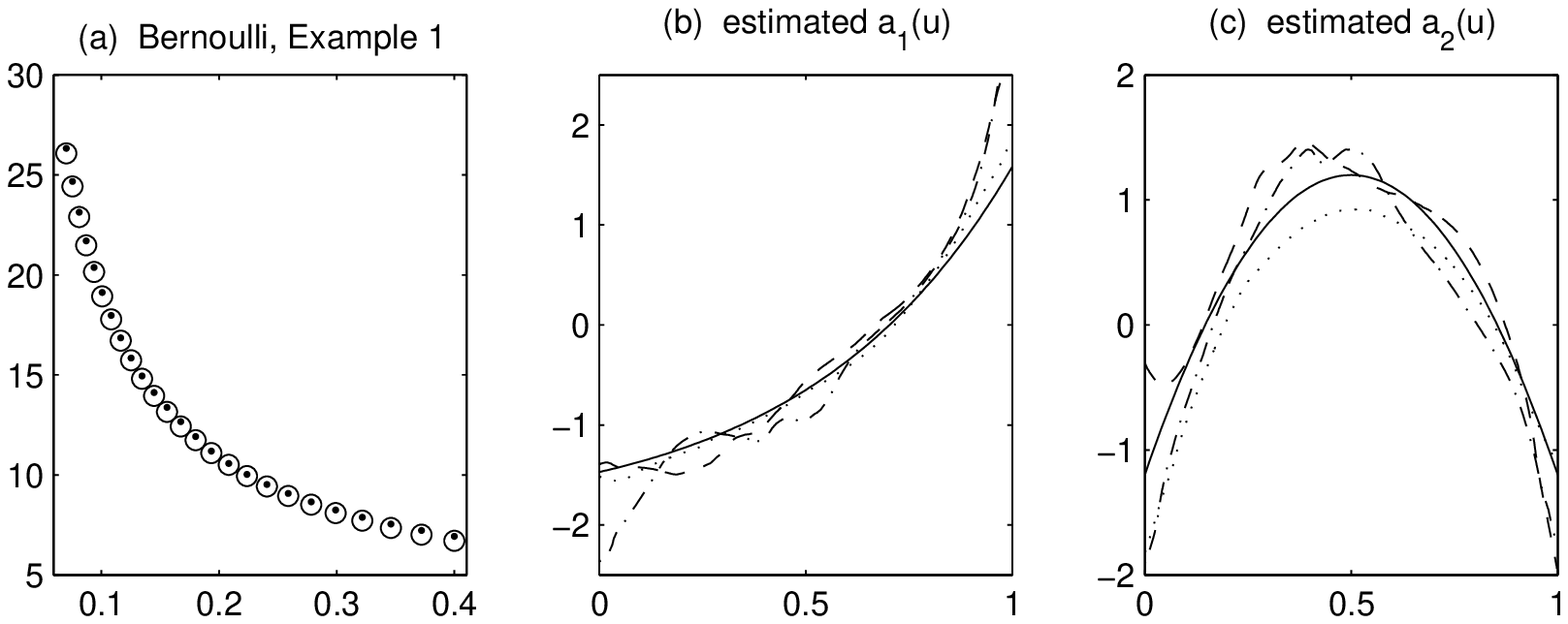,height=1.75in,width=6in}}
\vskip 0.1in
\centerline{\psfig{figure=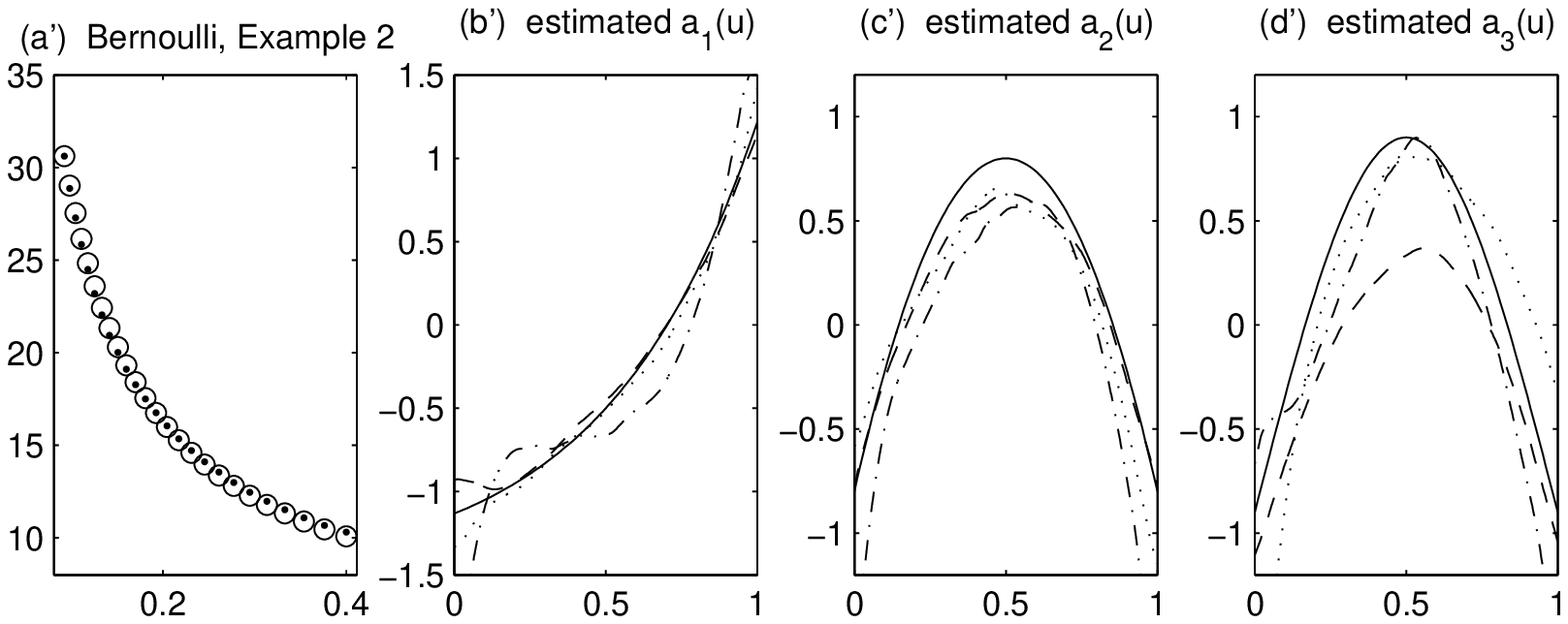,height=1.75in,width=6in}}
\begin{singlespace}
\caption {\textsl{Evaluation of Local-Likelihood Varying Coefficient
Regression for Bernoulli Responses. Captions are similar to those
for Figure \ref{figure-5}. Here $\hECV$ minimizes the empirical
version of (\ref{e5}); the empirical formula (\ref{e4}) uses
$\ta=.70$ and $\calC=1.09$ for $H_i^*$ and $\ta=.70$ and
$\calC=1.03$ for $\calS_i^*$.}} \label{figure-6}
\end{singlespace}
\end{figure}

{\bf Poisson regression:} We consider a variable $Y$, given
values $(u,\bvectorx)$ of the covariates $(U, \bvectorX)$,
following a Poisson distribution with parameter
$\lambda(u,\bvectorx)$, where the varying-coefficient functions
in $\ln\{\lambda(u,\bvectorx)\}$ are specified below:
\begin{eqnarray*}
\mbox{Example 1:} && d=2,\ a_1(u)=5.5+.1\, \exp(2u-1),\ a_2(u)=.8u(1-u), \\
\mbox{Example 2:} && d=3,\ a_1(u)=5.5+.1\, \exp(2u-1),\
a_2(u)=.8u(1-u),\ a_3(u)=.2\sin^2(2\pi u).
\end{eqnarray*}
We assume that $U$ is a uniform random variable on the interval
$[0,1]$ and is independent of $\bvectorX=(X_1, X_2, X_3)^T$, with
$X_1\equiv 1$, where $(X_2, X_3)$ follows a zero-mean and
unit-variance bivariate normal distribution with correlation
coefficient $1/\sqrt 2$. In Figure \ref{figure-5}, plot (a) reveals
that the actual degrees of freedom are well captured by the
empirical formula (\ref{e4}). To evaluate the performance of
$\hECV$, we generate $100$ random samples of size $400$. Figure
\ref{figure-5}(b)--(c) plot the estimated curves of $a_1(u)$ and
$a_2(u)$ from three typical samples. The typical samples are
selected so that their ASE values, in which
$\mbox{ASE}=n^{-1}\sum_{j=1}^{n} \{\hat
\theta(U_j,\bvectorX_j)-\theta(U_j,\bvectorX_j)\}^2$, correspond to
the $25$th (dotted line), $50$th (dashed line), and $75$th
(dash-dotted line) percentiles in the $100$ replications. The
corresponding results for Example 2 are given in Figure
\ref{figure-5}(a')--(d'). These plots provide convincing evidences
that $\hECV$, when applied to recovering multiple smooth curves
simultaneously, performs competitively well with that to fitting a
single smooth curve.
\begin{center}
{[ \textsl{Put Figure \ref{figure-6} about here} ]}
\end{center}

{\bf Logistic regression:} We now consider the varying-coefficient
logistic regression model for Bernoulli responses, in which
varying-coefficient functions in
$\logit\{P(Y=1|U=u,\bvectorX=\bvectorx)\}$ are specified as
\begin{eqnarray*}
\mbox{Example 1:} && d=2,\ a_1(u)=1.3\{\exp(2u-1)-1.5\},\ a_2(u)=1.2\{8u(1-u)-1\}, \\
\mbox{Example 2:} && d=3,\ a_1(u)=\exp(2u-1)-1.5,\
a_2(u)=.8\{8u(1-u)-1\},\ a_3(u)=.9\{2\sin(\pi u)-1\}.
\end{eqnarray*}
We assume that $X_1=1$; $X_2$ and $X_3$ are uncorrelated standard
normal variables, and are independent of $U\sim U(0,1)$. Figure
\ref{figure-6} depicts plots whose captions are similar to those
for Figure \ref{figure-5}. Compared with previous examples of
univariate logistic regression and varying-coefficient Poisson
regression, the current model fitting for binary responses is
considerably more challenging. Despite the increased difficulty,
the LB local-likelihood logistic regression estimates, through the
use of the hybrid bandwidth selector $\hECV$, captures the major
features of the model structure with reasonably good details.

\section{Real Data Applications} \label{sect-8}

In this section, we apply the hybrid bandwidth selection method
for binary responses to analyze an employee dataset (Example 11.3
of Albright, \etal 1999) of the Fifth National Bank of
Springfield, based on year $1995$ data.  The bank, whose name has
been changed, was charged in court with that its female employees
received substantially smaller salaries than its male employees.
For each of its $208$ employees, the dataset consists of eight
variables, including
\begin{itemize}
\item \textsf{JobGrade}: a categorical variable for the current job level,
      with possible values 1--6 (6 is highest)
\item \textsf{YrHired}: year  employee was hired
\item \textsf{YrBorn}: year  employee was born
\item \textsf{Gender}: a categorical variable with values ``Female" and ``Male"
\item \textsf{YrsPrior}: number of years of work experience at another
        bank prior to working at Fifth National.
\end{itemize}
The data set was carefully analyzed by Fan and Peng (2004).  After
adjusting for covariates such as age, years of work experience, and
education level, they did not find stark evidence of discrimination.
However, they pointed out that $77.29\%$ ($R^2$) of the salary
variation can be explained by the job grade alone and the question
becomes whether it was harder for females to be promoted, after
adjusting for confounding variables such as age and years of work
experience.  They did not carry out the analysis further.

\begin{figure}[htbp]
\centerline{\psfig{figure=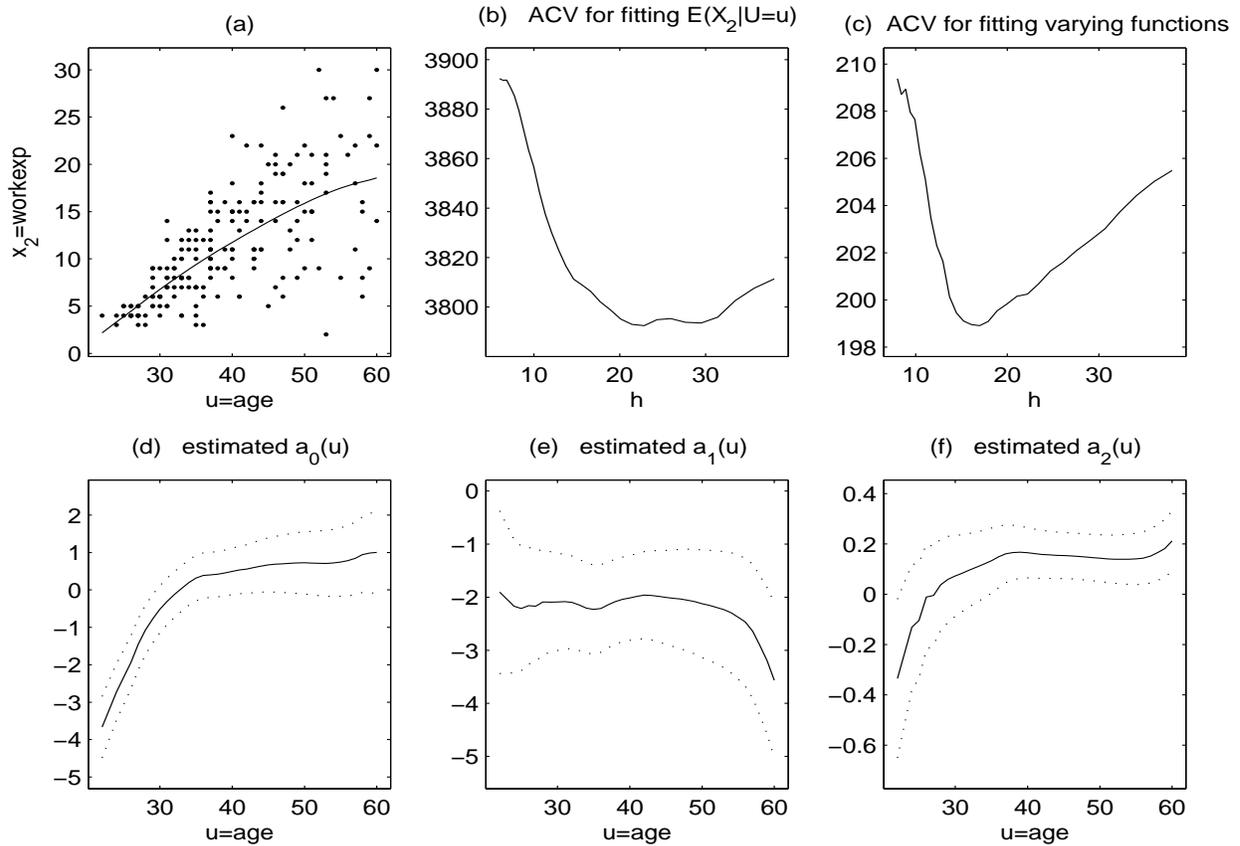,height=4.5in,width=6.5in}}
\begin{singlespace}
\caption {\textsl{Applications to the Job Grade Data Set Modeled by
(\ref{h1}). (a) scatter plot of work experience versus age along
with a local linear fit; (b) plot of the approximate CV function
against bandwidth for the local linear fit in (a); (c) plot of the
approximate CV function, defined in (\ref{e5}), against bandwidth
     for fitting varying coefficient functions;
(d) estimated $a_0(u)$; (e) estimated $a_1(u)$; (f) estimated
$a_2(u)$. The dotted curves are the estimated functions plus/minus
$1.96$ times of the estimated standard errors.}} \label{figure-7}
\end{singlespace}
\end{figure}

To understand how the probability of promotion to high levels of
managerial job (and thus high salary) is associated with gender
and years of work experience, and how this association changes
with respect to age, we fit a varying-coefficient logistic
regression model,
\begin{eqnarray} \label{h1}
\logit\{P(Y=1|U=u,X_1=x_1,X_2=x_2)\}=a_0(u)+a_1(u)x_1+a_2(u)x_2,
\end{eqnarray}
with $Y$ the indicator of \textsf{JobGrade} at least $4$, $U$ the
covariate \textsf{Age}, $X_1$ the indicator of being
\textsf{Female}, and $X_2$ the covariate \textsf{WorkExp}
(calculated as $95-\textsf{YrHired}+\textsf{YrsPrior}$). Following
Fan and Peng (2004), outliers have been deleted, with the remaining
$199$ data for analysis. For this medium-lengthed data, use of the
bandwidth selector $\hACV$ which minimizes (\ref{e5}) seems to be
more natural than $\hECV$.
\begin{center}
{[ \textsl{Put Figure \ref{figure-7} about here} ]}
\end{center}

Our preliminary study shows a monotone decreasing pattern in the
fitted curve of $a_2(u)$. This is no surprise; the covariates
\textsf{Age} and \textsf{WorkExp} are highly correlated, as can be
viewed from the scatter plot in Figure \ref{figure-7}(a). Such
high correlation may cause some identifiability problem, thus in
model (\ref{h1}), we replace $X_2$ with a de-correlated variable,
$X_2-E(X_2|U)$, which is known to be uncorrelated with any
measurable function of $U$. The projection part, $E(X_2|U=u)$, can
easily be estimated by a univariate local linear regression fit.
Likewise, its bandwidth parameter can simply be chosen to minimize
the approximate cross-validation function (for Gaussian family),
illustrated in Figure \ref{figure-7}(b).

\begin{figure}[htbp]
\centerline{\psfig{figure=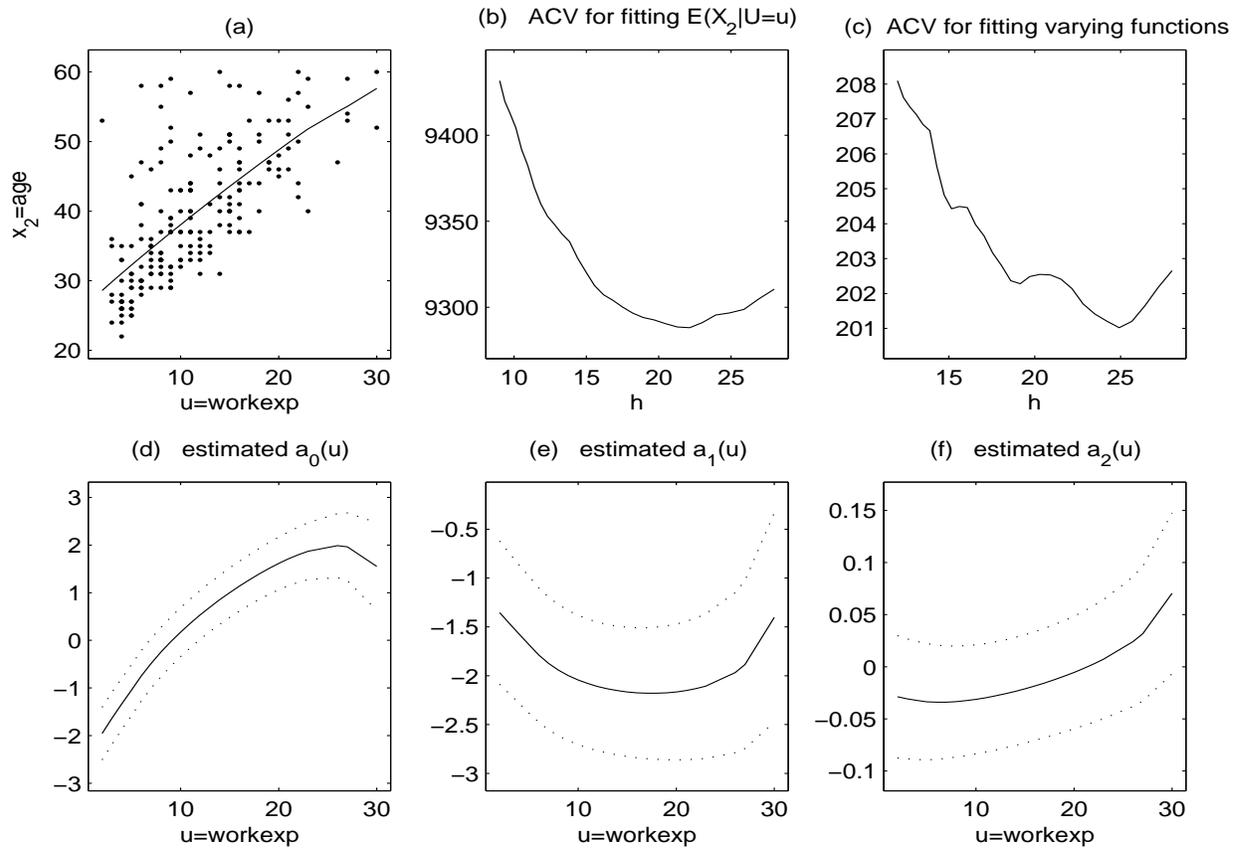,height=4.5in,width=6.5in}}
\begin{singlespace}
\caption {\textsl{Applications to the Job Grade Data Set Modeled by
(\ref{h1}). Captions are similar to those for Figure \ref{figure-7},
except that $U$ is \textsf{WorkExp} and $X_2$ is the de-correlated
\textsf{Age}.}} \label{figure-8}
\end{singlespace}
\end{figure}

After the de-correlation step, we now refit model (\ref{h1}). The
bottom panel of Figure \ref{figure-7} depicts the estimated varying
coefficient functions of $a_0(u)$, $a_1(u)$ and $a_2(u)$, plus/minus
the pointwise $1.96$ times of their estimated standard error. The
selected bandwidth is $16.9$ [see Figure \ref{figure-7}(c)]. Both
the intercept term and (de-correlated) \textsf{WorkExp} have the
statistically significant effects on the probability of promotion.
As an employee gets older, the probability of getting promoted keeps
increasing until around $40$ years of age and levels off after that.
It is interesting to note that the fitted coefficient function of
$a_1(u)$ for gender is below zero within the entire age span. This
may be interpreted as the evidence of discrimination against female
employees being promoted and lends support to the plaintiff.
\begin{center}
{[ \textsl{Put Figure \ref{figure-8} about here} ]}
\end{center}

To see whether the choice of smoothing variable $U$ makes a
difference in drawing the above conclusion, we fit again model
(\ref{h1}) with $U$ given by the covariate \textsf{WorkExp} and
$X_2$ by the de-correlated \textsf{Age} (due to the same reason of
monotonicity as in the previous analysis). Again, Figure
\ref{figure-8} shows that gender has an adverse effect and the
evidence for discrimination continues to be strong.  Indeed, the
estimated varying-function of $a_1(u)$ is qualitatively the same as
that in Figure \ref{figure-7}, as far as the evidence of
discrimination is concerned.

We would like to make a final remark on the de-correlation
procedure: This step does not alter (\ref{h1}),
particularly the function $a_1(\cdot)$. If this step is not taken,
then the estimate of $a_1(u)$ from either choice of $U$ continues to
be below zero and thus does not alter our previous interpretation
of the gender effect.

%\newpage
\begin{center}
\bf \Large References
\end{center}
\begin{description}
\item Albright, S. C., Winston, W. L., and Zappe, C. J. (1999),
      \textsl{Data Analysis and Decision Making with Microsoft Excel},
      Duxbury Press, Pacific Grove, California.

\item Altman, N., and MacGibbon, B. (1998),
      ``Consistent Bandwidth Selection for Kernel Binary Regression,''
      \textsl{J. Statist. Plann. Inference}, 70, 121--137.

\item Aragaki, A., and Altman, N. S. (1997),
      ``Local Polynomial Regression for Binary Response,''
      in \textsl{Computer Science and Statistics: Proceedings
      of the 29th Symposium on the Interface}.

\item Bickel, P.J. and Kwon, J. (2001),
      ``Inference for Semiparametric Models: Some Current Frontiers" (with discussion),
      \textsl{Statistica Sinica}, {11}, 863--960.

\item B\"ohning, D., and Lindsay, B. G. (1988),
      ``Monotonicity of Quadratic Approximation Algorithms,''
      \textsl{Ann. Inst. Statist. Math.}, 40, 641--663.

\item Br\`egman, L. M. (1967),
      ``A Relaxation Method of Finding a Common Point of Convex Sets and Its
      Application to the Solution of Problems in Convex Programming,"
      \textsl{U.S.S.R. Comput. Math. and Math. Phys.}, 7, 620--631.

\item Cai, Z., Fan, J. and Li, R. (2000),
     ``Efficient Estimation and Inferences for Varying-Coefficient
         Models," \textsl{Jour. Ameri. Statist. Assoc.}, { 95},
     888--902.

\item Efron, B. (2004),
      ``The Estimation of Prediction Error: Covariance Penalties and Cross-Validation"
      (with discussion),
      \textsl{J. Amer. Statist. Assoc.}, 99, 619--642.

\item Fan, J., and Chen, J. (1999),
      ``One-Step Local Quasi-Likelihood Estimation,''
      \textsl{J. R. Statist. Soc.}, Ser. B, 61, 927--943.

\item Fan, J., Heckman, N., and Wand, M. P. (1995),
      ``Local Polynomial Kernel Regression for Generalized
      Linear Models and Quasi-Likelihood Functions,''
      \textsl{J. Amer. Statist. Assoc.}, {90}, 141--150.

\item Fan, J., Farmen, M., and Gijbels, I. (1998),
      ``Local Maximum Likelihood Estimation and Inference,''
      \textsl{J. R. Statist. Soc.}, Ser. B, 60, 591--608.

\item Fan, J. and Huang, T. (2005),  ``Profile Likelihood Inferences
      on Semiparametric Varying-Coefficient Partially Linear
      Models,''
      \textsl{Bernoulli}, to appear.

\item Fan, J. and Peng, H. (2004),
      ``On Non-Concave Penalized Likelihood with Diverging Number of Parameters,"
      \textsl{The Annals of Statistics}, {32}, 928--961.

\item Golub, G. H., and Van Loan, C. F. (1996),
      \textsl{Matrix Computations} (Third edition),
      Baltimore, MD: Johns Hopkins University Press.

\item Hall, P. and Johnstone, I. (1992),
      ``Empirical Functionals and
      Efficient Smoothing Parameter Selection" (with discussion).
      \textsl{J. Royal. Statist. Soc. B}, {54}, 475--530.

\item Hardy, G. H., Littlewood, J. E., and P\'olya, G. (1988),
      \textsl{Inequalities} (Second edition),
      Cambridge, England: Cambridge University Press.

\item H\"ardle, W., Hall, P., and Marron, J. S. (1992),
      ``Regression Smoothing Parameters That Are Not Far From Their Optimum,"
      \textsl{J. Amer. Statist.  Assoc.}, {87}, 227--233.

\item Hastie, T. J., and Tibshirani, R. J. (1990),
      \textsl{Generalized Additive Models},
      London: Chapman and Hall.

\item Hastie, T. J., Tibshirani, R. and Friedman, J. (2001),
      \textsl{The Elements of Statistical Learning:  Data Mining,
      Inference, and Prediction}, Springer-Verlag, New York.

\item McCullagh, P., and Nelder, J. A. (1989),
      \textsl{Generalized Linear Models} (2nd ed.),
      London: Chapman and Hall.

\item Mitrinovi\'c, D. S., Pe\v cari\'c, J. E., and Fink, A. M. (1993),
      \textsl{Classical and New Inequalities in Analysis},
      Kluwer Academic Publishers Group, Dordrecht.

\item M\"uller, H.-G., and Schmitt, T. (1988),
      ``Kernel and Probit Estimates in Quantal Bioassay,''
      \textsl{J. Amer. Statist. Assoc.}, 83, 750--759.

\item Nelder, J. A., and Wedderburn, R. W. M. (1972),
      ``Generalized Linear Models,''
      \textsl{J. R. Statist. Soc.}, Ser. A, 135, 370--384.

\item Pregibon, D. (1981),
      ``Logistic Regression Diagnostics,''
      \textsl{Ann. Statist.}, 9, 705--24.

\item Rice, J. (1984),
      ``Bandwidth Choice for Nonparametric Regression,"
      \textsl{Ann. Statist.}, {20}, 712--736.

\item Ruppert, D. and Wand, M.P. (1994),
      ``Multivariate Weighted Least
      Squares Regression,"
      \textsl{Ann. Statist.}, {22}, 1346--1370.

\item Severini, T. A., and Staniswalis, J. G. (1994),
      ``Quasi-Likelihood Estimation in Semiparametric Models,''
      \textsl{J. Amer. Statist.  Assoc.}, 89, 501--511.

\item Shen, X., Huang, H.-C. and Ye, J. (2004),
      ``Adaptive Model Selection and Assessment for Exponential-Family Models,''
      \textsl{Technometrics}, 46, 306-317.

\item Staniswalis, J. G. (1989),
      ``The Kernel Estimate of a Regression Function in Likelihood-Based Models,''
      \textsl{J. Amer. Statist.  Assoc.}, 84, 276--283.

\item Tibshirani, R., and Hastie, T. (1987),
      ``Local Likelihood Estimation,''
      \textsl{J. Amer. Statist.  Assoc.}, 82, 559--567.

\item Tibshirani, R. (1996),
      ``Bias, Variance and Prediction Error for Classification Rules,''
      Technical report, Statistics Department, University of Toronto.

\item Wong, W. H. (1983),
      ``On the Consistency of Cross-Validation in Kernel Nonparametric Regression,"
      \textsl{Ann. Statist.}, 11, 1136--1141.

\item Yatchew, A. (1997),
      ``An Elementary Estimator for the Partial Linear Model,"
      \textsl{Economics Letters}, {57}, 135--143.

\item Zhang, C. M. (2003),
      ``Calibrating the Degrees of Freedom for Automatic Data Smoothing and
      Effective Curve Checking,''
      \textsl{J. Amer. Statist. Assoc.}, 98, 609--628.
\end{description}

\vskip 0.0in

\setcounter{equation}{0}
\renewcommand{\theequation} {A.\arabic{equation}}
\begin{center}
\bf \Large Appendix: Proofs of Main Results
\end{center}
We first impose some technical conditions. They are not the
weakest possible.

\noindent {\bf Condition (A)}:
\begin{itemize}
\item [(A1)] The function $q$ is concave and $q''(\cdot)$ is continuous.
\item [(A2)] $b''(\cdot)$ is continuous and bounded away from zero.
\item [(A3)] The kernel function $K$ is a symmetric probability density
       function with bounded support, and is Lipschitz continuous.
\end{itemize}

\noindent {\bf Condition (B)}:
\begin{itemize}
\item [(B1)] The design variable $X$ has a bounded support $\Omega_X$ and
      the density function $f_X$ which is Lipschitz continuous
      and bounded away from $0$.
\item [(B2)] $\theta(x)$ has the continuous $(p+1)$-th derivative in $\Omega_X$.
\end{itemize}

\noindent {\bf Condition (C)}:
\begin{itemize}
\item [(C1)] The covariate $U$ has a bounded support $\Omega_U$ and
        its density function $f_U$ is Lipschitz continuous
      and bounded away from $0$.
\item [(C2)] $a_j(u)$, $j=1,\ldots,d$, has the continuous $(p+1)$-th derivative in $\Omega_U$.
\item [(C3)] For the use of canonical links,
      the matrix $\Gamma(u)=E\{b''(\theta(u,\bvectorX))\bvectorX\bvectorX^T|U=u\}$ is positive
      definite for each $u\in \Omega_U$ and is Lipschitz continuous.
\end{itemize}

\noindent {\bf Notations}: Throughout our derivations, we
simplify notations by writing
$\theta_j(x;\bbeta)=\bx_j(x)^T\bbeta$,
$m_j(x;\bbeta)=b'(\theta_j(x;\bbeta))$,
$Z_j(x;\bbeta)=\{Y_j-m_j(x;\bbeta)\}/b''(\theta_j(x;\bbeta))$,
$\bz(x;\bbeta) = (Z_1(x;\bbeta),\ldots,Z_n(x;\bbeta))^T$, and
$w_j(x;\bbeta)=K_h(X_j-x) b''(\theta_j(x;\bbeta))$; their
corresponding quantities evaluated at $\hat \bbeta(x)$ are denote
by $\hat \theta_j(x)$, $\hat m_j(x)$, $\hat Z_j(x)$, $\hat
\bz(x)$, and $\hat w_j(x)$. Similarly, define $\hat
S_n(x)=S_n(x;\hat \bbeta(x))$.

Before proving the main results of the paper, we need the
following lemma.
\begin{lemma} \label{lem-2}
Define $\bV_i(\delta)=\diag\{\delta_{i1}, \ldots, \delta_{in}\}$.
For $\ell_{i,\delta} (\bbeta; x)$ defined in $(\ref{c2})$,
\begin{eqnarray}
\triangledown \ell_{i,\delta} (\bbeta; x)
&=& \bX(x)^T \bV_i(\delta) \bW(x;\bbeta) \bz(x;\bbeta)/a(\psi), \label{A1} \\
\triangledown^2 \ell_{i,\delta}(\bbeta; x) &=& -\bX(x)^T
\bV_i(\delta) \bW(x;\bbeta) \bX(x)/a(\psi), \label{A2}
\end{eqnarray}
in which
\begin{eqnarray}
\bX(x)^T \bV_i(\delta) \bW(x;\bbeta) \bz(x;\bbeta)
&=& \bX(x)^T \bW(x;\bbeta) \bz(x;\bbeta)-(1-\delta)\bx_i(x) K_h(X_i-x)\{Y_i-m_i(x;\bbeta)\},\qquad \label{A3} \\
\bX(x)^T \bV_i(\delta) \bW(x;\bbeta) \bX(x) &=& \bX(x)^T
\bW(x;\bbeta) \bX(x)-(1-\delta)w_i(x;\bbeta) \bx_i(x) \bx_i(x)^T.
\label{A4}
\end{eqnarray}
\end{lemma}
\textsl{Proof}: Defining a vector
$\btheta(x;\bbeta)=(\theta_1(x;\bbeta), \ldots,
\theta_n(x;\bbeta))^T$, we have that
\begin{eqnarray}
\triangledown \ell_{i,\delta}(\bbeta; x) &=& \frac {\partial
\btheta(x;\bbeta)} {\partial \bbeta}
    \frac {\partial \ell_{i,\delta}(\bbeta; x)}
    {\partial \btheta(x;\bbeta)}=\bX(x)^T \frac {\partial \ell_{i,\delta}(\bbeta; x)}
    {\partial \btheta(x;\bbeta)}, \label{A5}\\
\triangledown^2 \ell_{i,\delta}(\bbeta; x) &=& \bX(x)^T \frac
{\partial^2 \ell_{i,\delta}(\bbeta; x)}
    {\partial \btheta(x;\bbeta) \partial \btheta(x;\bbeta)^T} \bX(x). \label{A6}
\end{eqnarray}
Since
$$
\ell_{i,\delta}(\bbeta; x)=\sum_{j=1}^n \delta_{ij} [\{Y_j
\bx_j(x)^T\bbeta-b(\bx_j(x)^T\bbeta)\}/a(\psi)+c(Y_j,\psi)]K_h(X_j-x),
$$
it is easy to check that
\begin{eqnarray}
\frac {\partial \ell_{i,\delta}(\bbeta; x)} {\partial
\theta_j(x;\bbeta)} = \delta_{ij}\{Y_j-b'(\theta_j(x;\bbeta))\}
K_h(X_j-x)/a(\psi). \label{A7}
\end{eqnarray}
This combined with (\ref{A5}) leads to (\ref{A1}).

Following (\ref{A7}), we see that
\begin{eqnarray}
\frac {\partial^2 \ell_{i,\delta}(\bbeta; x)} {\partial
\theta_j(x;\bbeta) \partial \theta_k(x;\bbeta)} =0, \ j\ne k,
\quad \mbox{and}\quad \frac {\partial^2 \ell_{i,\delta}(\bbeta;
x)} {\{\partial \theta_j(x;\bbeta)\}^2} =-\delta_{ij}
b''(\theta_j(x;\bbeta)) K_h(X_j-x)/a(\psi). \label{A8}
\end{eqnarray}
This along with (\ref{A6}) indicates (\ref{A2}).

(\ref{A3}) and (\ref{A4}) can be obtained by decomposing the
identity matrix $\bI$ into $\bV_i(\delta)$ and
$\bI-\bV_i(\delta)$.

\subsection*{Proof of Proposition \ref{prop1}}

From (\ref{A1}) and (\ref{A2}), (\ref{c3}) can be rewritten as
\begin{eqnarray} \label{A9}
\bbeta_{L}=\bbeta_{L-1}+\{\bX(x)^T \bV_i(\delta)
\bW(x;\bbeta_{L-1}) \bX(x)\}^{-1} \{\bX(x)^T \bV_i(\delta)
\bW(x;\bbeta_{L-1}) \bz(x;\bbeta_{L-1})\}.
\end{eqnarray}
Setting $\delta=0$ in (\ref{A9}), the one-step estimate of $\hat
\bbeta^{-i}(x)$, which starts from $\bbeta_0=\hat \bbeta(x)$, is
given by
\begin{eqnarray} \label{A10}
\hat \bbeta(x)+ \{\bX(x)^T \bV_i(0) \bW(x;\hat \bbeta(x))
\bX(x)\}^{-1} \{\bX(x)^T \bV_i(0) \bW(x;\hat \bbeta(x)) \hat
\bz(x)\}.
\end{eqnarray}
Using the definition of $\hat \bbeta(x)$ (satisfying
$\triangledown \ell_{i,\delta}(\bbeta; x)=0$ with $\delta=1$),
along with (\ref{A3}) and (\ref{A4}), the above one-step estimate
of $\hat \bbeta^{-i}(x)$ equals
\begin{eqnarray}
\hat \bbeta(x)-\{\hat S_n(x)-\hat w_i(x)\bx_i(x)\bx_i(x)^T\}^{-1}
\bx_i(x) K_h(X_i-x) \{Y_i-\hat m_i(x)\}. \label{A11}
\end{eqnarray}
According to the Sherman-Morrison-Woodbury formula (Golub and Van
Loan 1996, p. 50),
\begin{eqnarray*}
\{\hat S_n(x)-\hat w_i(x)\bx_i(x)\bx_i(x)^T\}^{-1} = \{\hat
S_n(x)\}^{-1}
   +\frac {\hat w_i(x)\{\hat S_n(x)\}^{-1}
    \bx_i(x) \bx_i(x)^T \{\hat S_n(x)\}^{-1}}
    {1-\hat w_i(x) \bx_i(x)^T \{\hat S_n(x)\}^{-1} \bx_i(x)}.
\end{eqnarray*}
Thus
\begin{eqnarray*}
&& \{\hat S_n(x)-\hat w_i(x)\bx_i(x)\bx_i(x)^T\}^{-1} \bx_i(x)=\{\hat S_n(x)\}^{-1} \bx_i(x) \\
&& \qquad\qquad +\frac {\hat w_i(x)\{\hat S_n(x)\}^{-1} \bx_i(x)
\bx_i(x)^T
   \{\hat S_n(x)\}^{-1} \bx_i(x)}
   {1-\hat w_i(x) \bx_i(x)^T \{\hat S_n(x)\}^{-1} \bx_i(x)}
   =\frac {\{\hat S_n(x)\}^{-1} \bx_i(x)} {1- \calH_{ii}(x;\hat \bbeta(x))},
\end{eqnarray*}
by which (\ref{A11}) becomes
\begin{eqnarray*}
\hat \bbeta(x) -\frac {\{\hat S_n(x)\}^{-1} \bx_i(x) K_h(X_i-x)
\{Y_i-\hat m_i(x)\}}{1-\calH_{ii}(x;\hat \bbeta(x))}.
\end{eqnarray*}
This expression approximates $\hat \bbeta^{-i}(x)$ and thus leads
to (\ref{c6}).

Note that $\hat \theta_i=\hat \theta_i(X_i)$, $\hat m_i=\hat
m_i(X_i)$, $\hat \theta_i^{-i}=\hat \theta_i^{-i}(X_i)$, $\hat
m_i^{-i}=\hat m_i^{-i}(X_i)$, and
\begin{eqnarray} \label{A12}
H_i=\bevector_1^T\{\hat S_n(X_i)\}^{-1}\bevector_1 K_h(0)b''(\hat
\theta_i).
\end{eqnarray}
Applying (\ref{c6}), we have
\begin{eqnarray*}
& & \hat \theta^{-i}_i-\hat \theta_i
    =\bevector_1^T \{\hat \bbeta^{-i}(X_i)-\hat \bbeta(X_i)\} \\
& & \qquad \doteq -\frac {\bevector_1^T\{\hat S_n(X_i)\}^{-1}
\bevector_1 K_h(0) (Y_i-\hat m_i)}
    {1-\calH_{ii}(X_i;\hat \bbeta(X_i))}
    =-\frac {H_i} {1-H_i} (Y_i-\hat m_i)/b''(\hat \theta_i),
\end{eqnarray*}
leading to (\ref{c7}). This, together with a first-order Taylor's
expansion and the continuity of $b''$, yields
\begin{eqnarray*}
\hat m^{-i}_i-\hat m_i=b'(\hat \theta^{-i}_i)-b'(\hat \theta_i)
\doteq (\hat \theta^{-i}_i-\hat \theta_i)b''(\hat \theta_i),
\end{eqnarray*}
and thus (\ref{c8}).

\subsection*{Proof of Proposition \ref{prop2}}

By a first-order Taylor expansion, we have that
\begin{eqnarray*}
\hat \lambda_i-\hat \lambda_i^{-i}
&=& 2^{-1}\{q'(\hat m_i^{-i})-q'(\hat m_i)\}\doteq 2^{-1}q''(\hat m_i)(\hat m_i^{-i}-\hat m_i), \\
Q(\hat m_i^{-i},\hat m_i) &\doteq& -2^{-1}q''(\hat m_i)(\hat
m_i^{-i}-\hat m_i)^2.
\end{eqnarray*}
These, applied to an identity given in the Lemma of Efron (2004,
Section 4),
\begin{eqnarray*}
Q(Y_i,\hat m^{-i}_i)-Q(Y_i,\hat m_i)=2(\hat \lambda_i-\hat
\lambda_i^{-i}) (Y_i-\hat m_i^{-i})-Q(\hat m_i^{-i},\hat m_i),
\end{eqnarray*}
lead to
\begin{eqnarray*}
Q(Y_i,\hat m^{-i}_i)-Q(Y_i,\hat m_i)
&\doteq& q''(\hat m_i)(\hat m_i^{-i}-\hat m_i)(Y_i-\hat m_i^{-i})+2^{-1}q''(\hat m_i)(\hat m_i^{-i}-\hat m_i)^2  \\
&=& 2^{-1}q''(\hat m_i) \{(Y_i-\hat m_i)^2-(Y_i-\hat m_i^{-i})^2\}.
\end{eqnarray*}
Summing over $i$ and using (\ref{c8}) and (\ref{b7}), we complete
the proof.

\subsection*{Proof of Proposition \ref{prop3}}

From (\ref{c12}) and (\ref{c13}), we see that
\begin{eqnarray} \label{A13}
\frac{\hAMPEC(q_2)} {\hAMISE} =\bigg[\frac{|\Omega_X|
\int_{\Omega_X} F(x)G(x) dx} {\int_{\Omega_X} F(x) dx
\int_{\Omega_X} G(x) dx}\bigg]^{1/(2p+3)}.
\end{eqnarray}

For part (a), it suffices to consider oppositely ordered $F$ and
$G$. In this case, by the Tchebychef's inequality (Hardy, Littlewood,
and P\'olya, 1988, p. 43 and 168), we obtain
$$
|\Omega_X| \int_{\Omega_X} F(x)G(x) dx \le \int_{\Omega_X} F(x)
dx \int_{\Omega_X} G(x) dx.
$$
Since $F\ge 0$ and $G\ge 0$, it follows that
$$
\frac{|\Omega_X|\int_{\Omega_X} F(x)G(x) dx}{\int_{\Omega_X} F(x)
dx \int_{\Omega_X} G(x) dx} \le 1,
$$
which along with (\ref{A13}) indicates that $\hAMPEC(q_2)\le
\hAMISE$.

To verify part (b), it can be seen that under its assumptions,
for a constant $C>0$, $F(x)=C/|\Omega_X|b''(\theta(x))$ is
oppositely ordered with $G(x)=\{b''(\theta(x))\}^{-1}$, and thus
the conclusion of part (a) immediately indicates the upper bound
$1$. To show the lower bound, we first observe that (\ref{A13})
becomes
\begin{eqnarray} \label{A14}
\frac{\hAMPEC(q_2)}{\hAMISE} =\bigg[\frac{|\Omega_X|^2}
{\int_{\Omega_X} b''(\theta(x))dx \int_{\Omega_X}
\{b''(\theta(x))\}^{-1} dx}\bigg]^{1/(2p+3)}.
\end{eqnarray}
Incorporating the Gr\"uss integral inequality (Mitrinovi\'c, Pe\v
cari\'c, and Fink 1993),
$$
\bigg|\frac 1 {|\Omega_X|}\int_{\Omega_X} F(x)G(x)
dx-\frac{1}{|\Omega_X|^2} \int_{\Omega_X} F(x)dx \int_{\Omega_X}
G(x)dx\bigg| \le \frac 1 4 (M_F-m_F)(M_G-m_G),
$$
where $M_F=\max_{x\in \Omega_X} F(x)$, $m_F=\min_{x\in \Omega_X}
F(x)$, $M_G=\max_{x\in \Omega_X} G(x)$, and $m_G=\min_{x\in
\Omega_X} G(x)$, we deduce that
$$
\int_{\Omega_X} b''(\theta(x))dx \int_{\Omega_X}
\{b''(\theta(x))\}^{-1} dx \le
\frac{(m_{b''}+M_{b''})^2}{4m_{b''}M_{b''}}|\Omega_X|^2.
$$
This applied to (\ref{A14}) gives the lower bound.

\subsection*{Proof of Proposition \ref{prop4}}

Define $\bH=\diag\{1,h,\ldots,h^p\}$. From (\ref{A12}), we have
\begin{eqnarray}
H_i
&=& \bevector_1^T\{\hat S_n(X_i)/b''(\hat \theta(X_i))\}^{-1}\bevector_1 K_h(0) \nonumber \\
&=& n^{-1} \bevector_1^T\bH^{-1}\{n^{-1}\bH^{-1}\hat
S_n(X_i)/b''(\hat \theta(X_i))\bH^{-1}\}^{-1}
    \bH^{-1}\bevector_1 K_h(0) \nonumber \\
&=& (nh)^{-1} \bevector_1^T \{n^{-1}\bH^{-1}\hat
S_n(X_i)/b''(\hat \theta(X_i))\bH^{-1}\}^{-1}
    \bevector_1 K(0), \label{A15}
\end{eqnarray}
where $\hat S_n(x)=\sum_{j=1}^{n} \bx_j(x) \bx_j(x)^T K_h(X_j-x)
b''(\bx_j(x)^T\hat \bbeta(x))$. By Taylor's expansion and the
continuity assumptions on $b''$ and $f_X$, it follows that for
$x\in \Omega_X$,
$$
n^{-1}\bH^{-1}\hat S_n(x)/b''(\hat
\theta(x))\bH^{-1}=f_X(x)S+o_P(1).
$$
Combining this expression with (\ref{A15}), it can be shown that
$$
\sum_{i=1}^n H_i =\sum_{i=1}^n \frac 1 {nh f_X(X_i)}
\bevector_1^T S^{-1} \bevector_1 K(0)\{1+o_P(1)\} =\frac
{\calK(0)} {nh} \sum_{i=1}^n \frac 1{f_X(X_i)} \{1+o_P(1)\},
$$
which will finish the proof.
\begin{lemma} \label{lem-3}
Assume that the kernel function $K$ is non-negative, symmetric
and uni-modal. Then for $i=1,\ldots,n$, $\calS_i$ is a decreasing
function of $h>0$ for which $\calS_i$ is well-defined.
\end{lemma}
\textsl{Proof}: Consider the matrices $A_i(h)=\bX(X_i)^T
\diag\{K(|X_j-X_i|/h)\}_{j=1}^{n} \bX(X_i)$, $i=1,\ldots,n$.
If $K$ is non-negative and uni-modal, then $0<h_1<h_2$ implies that
$A_i(h_1)\le A_i(h_2)$ or, equivalently, $\{A_i(h_1)\}^{-1} \ge
\{A_i(h_2)\}^{-1}$. We complete the proof by noting
$\calS_i=\bevector_1^T \{A_i(h)\}^{-1} \bevector_1 K(0)$, since
$K$ is symmetric.

\subsection*{Proof of Proposition \ref{prop5}}

The one-step estimate of $\hat \bbeta^{-i}(x)$, starting from
$\bbeta_0=\hat \bbeta(x)$, is given by
\begin{eqnarray} \label{A16}
\hat \bbeta(x)+ 4\{\bX(x)^T \bV_i(0) \bK(x) \bX(x)\}^{-1}
\{\bX(x)^T \bV_i(0) \bW(x;\hat \bbeta(x)) \hat \bz(x)\},
\end{eqnarray}
i.e., $\hat
\bbeta(x)-4\{S_n(x)-K_h(X_i-x)\bx_i(x)\bx_i(x)^T\}^{-1}\bx_i(x)K_h(X_i-x)\{Y_i-\hat
m_i(x)\}$. Again, using the Sherman-Morrison-Woodbury formula
(Golub and Van Loan 1996, p. 50),
\begin{eqnarray*}
\{S_n(x)-K_h(X_i-x)\bx_i(x)\bx_i(x)^T\}^{-1} = \{S_n(x)\}^{-1}
    +\frac {K_h(X_i-x)\{S_n(x)\}^{-1} \bx_i(x)\bx_i(x)^T \{S_n(x)\}^{-1}}
    {1-K_h(X_i-x) \bx_i(x)^T \{S_n(x)\}^{-1} \bx_i(x)},
\end{eqnarray*}
and thus
\begin{eqnarray*}
&& \{S_n(x)-K_h(X_i-x)\bx_i(x)\bx_i(x)^T\}^{-1}\bx_i(x)=\{S_n(x)\}^{-1} \bx_i(x) \\
&& \qquad +\frac {K_h(X_i-x)\{S_n(x)\}^{-1} \bx_i(x) \bx_i(x)^T
   \{S_n(x)\}^{-1} \bx_i(x)} {1-K_h(X_i-x) \bx_i(x)^T \{S_n(x)\}^{-1} \bx_i(x)}
   =\frac {\{S_n(x)\}^{-1} \bx_i(x)} {1- K_h(X_i-x) \bx_i(x)^T \{S_n(x)\}^{-1} \bx_i(x)},
\end{eqnarray*}
by which (\ref{A16}) becomes
\begin{eqnarray*}
\hat \bbeta(x) -\frac {4\{S_n(x)\}^{-1} \bx_i(x) K_h(X_i-x)
\{Y_i-\hat m_i(x)\}} {1-K_h(X_i-x) \bx_i(x)^T \{S_n(x)\}^{-1}
\bx_i(x)}.
\end{eqnarray*}
This expression approximates $\hat \bbeta^{-i}(x)$ and thus leads to
(\ref{d3}).

Applying (\ref{d3}), we have
\begin{eqnarray*}
& & \hat \theta^{-i}_i-\hat \theta_i
    =\bevector_1^T \{\hat \bbeta^{-i}(X_i)-\hat \bbeta(X_i)\} \\
& & \qquad \doteq -\frac {4\bevector_1^T\{S_n(X_i)\}^{-1}
\bevector_1 K_h(0) (Y_i-\hat m_i)}
    {1-\calS_i}
    =-\frac {4\calS_i} {1-\calS_i} (Y_i-\hat m_i),
\end{eqnarray*}
leading to (\ref{d4}). Proofs of (\ref{d5}) and (\ref{d6}) are
similar to those of Proposition \ref{prop2}.

\subsection*{Proofs of Propositions \ref{prop6}--\ref{prop7}}
The technical arguments are similar to the proofs of Propositions
\ref{prop1}--\ref{prop2} and thus details are omitted.

\subsection*{Proof of Proposition \ref{prop8}}
Recalling the definition of $H_i^*$ in Section \ref{sect-5.2}, we
have that
\begin{eqnarray}
H_i^* &=& (\bevector_1\otimes \bvectorX_i)^T
\{S_n^*(U_i;\hat\bbeta(U_i))\}^{-1} (\bevector_1\otimes
\bvectorX_i)\
    \{K_h(0)b''(\hat \theta(U_i,\bvectorX_i))\} \nonumber \\
&=& (nh)^{-1} (\bevector_1\otimes \bvectorX_i)^T
    \{n^{-1}(\bH\otimes \bI_d)^{-1}\hat S_n^*(U_i)(\bH\otimes \bI_d)^{-1}\}^{-1}
    (\bevector_1\otimes \bvectorX_i) K(0)b''(\hat \theta(U_i,\bvectorX_i)), \qquad \label{A17}
\end{eqnarray}
where
\begin{eqnarray*}
\hat S_n^*(u) &=& \sum_{j=1}^{n} \{\bu_j(u)\otimes
\bvectorX_j\}\{\bu_j(u)\otimes \bvectorX_j\}^T
    K_h(U_j-u) b''(\{\bu_j(u)\otimes \bvectorX_j\}^T\hat \bbeta(u)) \\
&=& \sum_{j=1}^{n} \big[\{\bu_j(u)\bu_j(u)^T\} \otimes
(\bvectorX_j\bvectorX_j^T)\big]
    K_h(U_j-u) b''(\{\bu_j(u)\otimes \bvectorX_j\}^T\hat \bbeta(u)).
\end{eqnarray*}
It can be shown that for $u\in \Omega_U$,
\begin{eqnarray*}
& & n^{-1} (\bH\otimes \bI_d)^{-1}\hat S_n^*(u)(\bH\otimes \bI_d)^{-1} \\
&=& n^{-1}
    \sum_{j=1}^{n}
    \big[\big\{\bH^{-1}\bu_j(u)\bu_j(u)^T\bH^{-1}\big\}\otimes (\bvectorX_j\bvectorX_j^T)\big]
    K_h(U_j-u) b''(\{\bu_j(u)\otimes \bvectorX_j\}^T\hat \bbeta(u)) \\
&=& n^{-1}
    \sum_{j=1}^{n}
    \big[\big\{\bH^{-1}\bu_j(u)\bu_j(u)^T\bH^{-1}\big\}\otimes (\bvectorX_j\bvectorX_j^T)\big]
    K_h(U_j-u) b''(\theta(u,\bvectorX_j))+o_P(1) \\
&=& f_U(u)[S \otimes
E\{b''(\theta(u,\bvectorX))\bvectorX\bvectorX^T|U=u\}]+o_P(1)
    =f_U(u)\{S \otimes \Gamma(u)\}+o_P(1).
\end{eqnarray*}
This expression applied to (\ref{A17}) further implies that
\begin{eqnarray}
\sum_{i=1}^n H_i^* &=& \sum_{i=1}^n \frac {1}{nhf_U(U_i)}
(\bevector_1\otimes \bvectorX_i)^T \{S^{-1} \otimes
\Gamma(U_i)^{-1}\}
    (\bevector_1\otimes \bvectorX_i) K(0) b''(\hat \theta(U_i,\bvectorX_i))\{1+o_P(1)\} \nonumber \\
&=& \sum_{i=1}^n \frac {1}{nhf_U(U_i)} (\bevector_1^T
S^{-1}\bevector_1)K(0)
    \big\{\bvectorX_i^T \Gamma(U_i)^{-1}\bvectorX_i b''(\hat \theta(U_i,\bvectorX_i))\big\} \{1+o_P(1)\} \nonumber \\
&=& \sum_{i=1}^n \frac {\calK(0)}{nhf_U(U_i)}
    \big\{\bvectorX_i^T \Gamma(U_i)^{-1}\bvectorX_i b''(\theta(U_i,\bvectorX_i))\big\}\{1+o_P(1)\}. \label{A18}
\end{eqnarray}
For (\ref{A18}), a direct calculation gives that
\begin{eqnarray*}
E\big\{\bvectorX^T \Gamma(U)^{-1}\bvectorX b''(\theta(U,\bvectorX))/f_U(U) \big\}
&=& \tr \big[E \big\{\Gamma(U)^{-1}\bvectorX \bvectorX^T b''(\theta(U,\bvectorX))/f_U(U) \big\}\big] \\
&=& \tr\big\{ E\big[\Gamma(U)^{-1}E\{b''(\theta(U,\bvectorX)) \bvectorX \bvectorX^T|U\}/f_U(U)\big] \big\} \\
&=& \tr\big[E\big\{\Gamma(U)^{-1}\Gamma(U)/f_U(U)\big\}\big] \\
&=& \tr(|\Omega_U|\bI_d)=d|\Omega_U|.
\end{eqnarray*}
This completes the proof.
\newpage
%-------------------------------------------------------------------------------------------------
% band_optimal_grm.m
\begin{table}[bhtp]
\centering
\begin{singlespace}
\caption{\textsl{The Asymptotic Optimal Bandwidths $\hAMPEC(q_2)$
Calculated From (\ref{c12}) and $\hAMISE$ From (\ref{c13}), with
$p=1$, the Epanechnikov Kernel, and Examples Given in Section
\ref{sect-7.1}}} \label{table-1}
\end{singlespace}
\begin{center}
\begin{tabular}{cccc} \hline \hline
\textsl{Exponential family} & \textsl{Example} & $\hAMPEC(q_2)$ &
$\hAMISE$ \\ \hline
Poisson   & 1  & $.070$ & $.079$  \\
          & 2  & $.089$ & $.099$  \\
          & 3  & $.127$ & $.136$  \\ \hline
Bernoulli & 1  & $.106$ & $.108$  \\
          & 2  & $.151$ & $.146$  \\
          & 3  & $.184$ & $.188$  \\ \hline
\end{tabular}
\end{center}
\end{table}
%
%------------------------------------------------------------------------------------------------
% choice_a.m
\begin{table}[bhtp]
\centering
\begin{singlespace}
\caption{\textsl{Choices of $\ta$ and $\calC$, in the Empirical
Formulas (\ref{c14}) and (\ref{e4}), for the $p$th Degree Local
Polynomial Regression for Gaussian Responses}} \label{table-2}
\end{singlespace}
\begin{center}
\begin{tabular}{cccc|cccc} \hline \hline
\textsl{Design type} & $p$ & $\ta$ & $\calC$ & \textsl{Design type}
& $p$ & $\ta$ & $\calC$ \\ \hline
        fixed   & $0$ & $0.55$ &   $1$   &      random  & $0$ & $0.30$ & $0.99$ \\
                & $1$ & $0.55$ &   $1$   &              & $1$ & $0.70$ & $1.03$ \\
                & $2$ & $1.55$ &   $1$   &              & $2$ & $1.30$ & $0.99$ \\
                & $3$ & $1.55$ &   $1$   &              & $3$ & $1.70$ & $1.03$ \\ \hline
\end{tabular}
\end{center}
\end{table}
%
%-------------------------------------------------------------------------------------------------
% margin.m
%
%----------------------------------------------------------------------------------
% Q_plot.m
%
%-------------------------------------------------------------------------------------------------
% ENP_LL_grm.m
%
%-----------------------------------------------
% ENP_LL_grm.m
%
%-----------------------------------------------
% ENP_LL_gvcm.m
%
%-----------------------------------------------
% ENP_LL_gvcm.m
%
%-----------------------------------------------
% bank.m
%
%-----------------------------------------------
% bank.m

\end{document}